\documentclass[11pt, a4paper, leqno]{article}
\usepackage{amssymb}
\usepackage{amscd}
\usepackage{amsmath}
\usepackage{amsthm}
\usepackage{mathrsfs}
\usepackage{graphics}
\usepackage[pdftex]{graphicx}
\newtheorem{rmk}[equation]{\indent \rm {\it Remark}}

\newcommand{\C}{{\mathbf{C}}}

\newcommand{\del}{{\partial}}

\renewcommand{\O}{{\mathcal{O}}}

\renewcommand{\P}{{\mathbf{P}}}

\newcommand{\PD}{{\mathrm{P}\Delta}}
\newcommand{\pnc}{{\mathbf{P}^n(\mathbf{C})}}

\newcommand{\ponec}{\P^1(\C)}

\newcommand{\R}{{\mathbf{R}}}

\newcommand{\Sig}{\mathit{\Sigma}}

\numberwithin{equation}{section}
\def\labelenumi{\rm(\roman{enumi})}

\title{\bf
On Kiyoshi Oka's Unpublished Papers in 1943
}
\date{For the 120th anniversary of Kiyoshi Oka's birth}
\author{
By J{\sc unjiro} N{\sc oguchi}\thanks{Research supported in part by Grant-in-Aid
 for Scientific Research (C) 19K03511. \hfill\break
\hbox{\quad} MSC2020: 32A99; 32E30; 01A60 \hfill\break
\hbox{\quad} Key words: coherence, Oka, Levi problem,  Hartogs' inverse problem,
 several complex variables
 \hfill\break
\hbox{\quad} Affiliation and address: {Graduate School of Mathematical Sciences},
{University of Tokyo (Emeritus)};
{Komaba, Meguro-ku, Tokyo 153-8914}, {Japan}\hfill\break
\hbox{\quad} e-mail: noguchi@ms.u-tokyo.ac.jp
}}
\begin{document}
%
\setlength{\baselineskip}{15pt}
\parskip+2.5pt
\maketitle
\thispagestyle{empty}

\begin{abstract}
In 1943 from September to December Kiyoshi Oka wrote a series
of papers numbered from VII to XI, as the research reports
to Teiji Takagi (then, Professor of Tokyo Imperial University),
in which he solved affirmatively the so-called Levi Problem
(Hartogs' Inverse Problem termed by Oka) for unramified
Riemann domains over $\C^n$. This problem which had been left open
for more than thirty years then, was the last one of the Three Big
 Problems summarized by Behnke--Thullen 1934.
 The papers were hand-written in Japanese,
 consist of pp.~108 in total, and have not been published
by themselves.
The aim of the present article is to provide an English translation
of the most important, last paper (Part II) with
 preparation (Part I).
At the end of Part I we will discuss a problem which K. Oka left
and is still open.
\end{abstract}

\bigskip
\centerline{\Large\bf  Part I}

\medskip
In this Part I we discuss Kiyoshi Oka's unpublished series of
five papers, VII---XI in 1943 (\cite{okap}),
which were hand-written in Japanese and consist of pp.~108 in total.
In Part II we present the English translation of the last
one XI of \cite{okap} that contains the most important main results.
{Part II is the main part of the present article.
In Part I it is not intended to survey the developments of the
subject since the time of Oka and thereafter, but rather is aimed to
serve for the preparations or a sort of appendices,
 so that Part II is readable for general readers
 without specific knowledge of the subject
at the time. Moreover, if one gets into the proofs described in XI,
he will still find methods that have not been presented in published references,
 so far by the author's knowledge, and are original and have interests
 even from the present
viewpoint.}
For general references 
{about the developments of the present subject},
 cf., e.g.,  Hitotsumatsu \cite{hi},
 Gunning--Rossi \cite{guro}, H\"ormander \cite{ho2},
Nishino \cite{nis96}, Lieb \cite{li},
 Noguchi \cite{nog16}, \cite{nog19}.

The method of the proof of the Pseudoconvexity Problem (i.e., Hartogs' Inverse
 Problem, Levis's Problem) given in this series of papers 1943
is quite similar to that of  Oka IX published in 1953
 except for the use of Coherence Theorems:
There, in the unpublished papers 1943,
 he proved some ideal theoretic properties of holomorphic
functions, which was sufficient to prove
 the J\^oku-Ik\^o (lifting principle)
 with estimates; then it led to the solution of
the Pseudoconvexity Problem.
 In this series of papers, he already had in mind a project
not only to settle the Pseudoconvexity Problem
of  general dimension, but also to deal with the problem
for ramified Riemann domains; it would actually lead to the notion
of ``Coherence''.

Reading the series of unpublished papers 1943
 we see the dawn of the
then unknown notion of ``Coherence'' or 
{\em ``Id\'eaux de domaines ind\'etermin\'es}''
in Oka's terms, and may observe that
the turn of years ``{\em 1943/'44}\,'' was indeed the watershed in
the study of analytic function theory of several variables. 

\section{Three Big Problems}\label{3bp}
\textbf{a)}
K. Oka's research \cite{oka}, I---IX (published) was motivated
by the monograph of Behnke--Thullen \cite{bt} 1934:
They summarized the main problems then in the theory of complex
analytic functions of several variables, listing the following
{\em Three Big Problems}.
\begin{enumerate}
\setlength{\itemsep}{-2pt}
\item
The Levi (Hartogs' Inverse) Problem.
\item
Cousin (I/II) Problem.
\item
Problem of expansions of functions (Approximation Problem).
\end{enumerate}

These problems are well-known among complex analysists, but
 we will recall  for convenience
 the above problems in the next subsection b),
following Behnke--Thullen \cite{bt} (cf.\ Lieb \cite{li}).

The difficulty of the problems was referred by H. Cartan \cite{oka2}
as ``{\em quasi-surhumaine (quasi-superhuman)}'' and by R. Remmert \cite{oka2}
as ``{\em Er  l\"oste Probleme, die als  unangreitbar galten
(He solved
problems which were believed to be unsolvable)}''.

K. Oka solved all these problems in the opposite order.
By establishing 
``J\^oku-Ik\^o''\footnote{\label{jkik}\,This consists of two
 (Japanese) words, and means that ``{\em one transfers himself from
the original space of the given dimension to a space of even higher
dimension''.} Cf.\ \S\ref{jokuiko}}
 in \cite{oka} I--II, he proved
Problem (iii) above and  (ii) the Cousin I Problem,
 and then in \cite{oka} III,
he obtained the Oka Principle, settling (ii) the Cousin II Problem.
The most difficult problem (i) was first proved for
univalent domains (subdomains) of $\C^2$ in \cite{oka} VI 1942, leaving
for the general dimensional case the last paragraph of the paper:
\vspace{-8pt}
\begin{quote}\sf
``{\em L'auteur pense que cette conclusion sera aussi ind\'ependante des  
nombres de variables complexes. (The author thinks that this conclusion
will be also independent of the number of complex variables.)}''
\end{quote}

But, it was a general cognition that the higher dimensional case
was still open (in Japan there seems to have been a sentiment that
the higher dimensional case of univalent domains was  already settled),
 and it was proved as follows:
{\def\labelenumi{\rm(\arabic{enumi})}
\begin{enumerate}
\item
S. Hitotsumatsu \cite{hi} (a short note in Japanese was published),
 1949 for univalent domains
of $\C^n~(n \geq 2$, same as in (iii) below by Weil's integral).
\item
K. Oka \cite{oka} IX, 1953 for unramified Riemann domains over
$\C^n$ (by Coherence, J\^oku-Ik\^o and Cauchy integral).
\item
H.J. Bremermann \cite{br} and F. Norguet \cite{nor} 1954,
 independently for univalent domains of $\C^n$ (by Weil's integral).
\end{enumerate} }

\textbf{b)} (i)
To get the idea of the problems we consider a univalent
domain (i.e., a subdomain) $\Omega$ of $\C^n$.
Let $\Omega' \supset \Omega$ be a domain of $\C^n$. 
If every holomorphic function in $\Omega$ is extendable
 to a holomorphic function
in $\Omega'$, $\Omega'$ is called an {\em extension of holomorphy}
of $\Omega$. In the case of $n=1$, there is no extension
of holomorphy other than $\Omega'=\Omega$, but in the case
case of $n \geq 2$, $\Omega' \supsetneq \Omega$ can happen
(Hartogs' phenomenon, 1906--).
For example, let  $n \geq 2$, 
let $a=(a_1, \ldots, a_n) \in \C^n$ and define
$\Omega_{\mathrm{H}}(a;\delta, \gamma) \subset \C^n$, so-called
a Hartogs domain, as follows:
With a pair of $n$-tuples of positive numbers,
$\gamma=(\gamma_j)_{1 \leq j \leq n}$ and
$\delta=(\delta_j)_{1 \leq j \leq n}$ satisfying
 $0<\delta_j< \gamma_j ~(1 \leq j \leq n)$, we set
\begin{align}
\label{hartogs}
\PD(a; \gamma) &=
 \{z=(z_1,\ldots, z_n) \in \C^n: |z_j-a_j| <\gamma_j, 1 \leq j
 \leq n\},\\
\notag
\Omega_1 &=\{z=(z_1,\ldots, z_n) \in \PD(a; \gamma):
  |z_j-a_j|<\delta_j,\: 2 \leq j \leq n\},\\
\nonumber
\Omega_2 &=\{z=(z_1,\ldots, z_n)\in \PD(a; \gamma):
  \delta_1 < |z_1-a_1| <\gamma_1\}, \\
\notag
\Omega_{\mathrm{H}}(a; \delta, \gamma) &= \Omega_1 \cup \Omega_2
 \subsetneq \PD(a; \gamma) .
\end{align}
It is immediate to see that the polydisk $\PD(a; \gamma)$ is an extension
of holomorphy of $\Omega_{\mathrm{H}}(a;\delta, \gamma)$
(cf., e.g., \cite{nog16} \S1.2.4).

The notion of the ``extension of holomorphy'' is naturally generalized
to the case of multi-sheeted (ramified or unramified) domains
over $\C^n$ and this is definitely necessary in the case of $n \geq 2$;
in fact, it is known that there is a subdomain of $\C^2$
which has an infinitely-sheeted unramified domain over $\C^2$
as an extension of holomorphy (cf., e.g., \cite{nog16} \S5.1).
In this paper, domains over $\C^n$ are {\em unramified},
as far as it is not mentioned to be ramified.

Now, let $\Omega$ be a domain over $\C^n$.
The maximal domain among the extensions of holomorphy of $\Omega$
is called the {\em envelope of holomorphy} of $\Omega$,
 denoted by $\hat\Omega$.
It exists, but is not necessarily univalent even if $\Omega$
is univalent as mentioned above.

If $\Omega=\hat\Omega$,
$\Omega$ is called a {\em domain of holomorphy}.
In the above example, $\PD(a; \gamma)$ is the envelope of holomorphy
of $\Omega_{\mathrm{H}}(a;\delta, \gamma)$ and a domain of holomorphy.
 Hartogs' phenomenon implies that the shape of
 singularities of holomorphic functions
is not arbitrary; contrarily, before Hartogs it had been thought
arbitrary. In the study of the shape of singularities of holomorphic
functions, in other words, the shape of the boundary of
 a domain of holomorphy $\Omega$,
E.E. Levi found around 1910 in the case of $n=2$
 that with the assumption of the $C^2$-regularity of
the boundary $\del \Omega$ defined by $\varphi$ so that
$\Omega=\{\varphi <0\}$, $d\varphi \not=0$ on $\del \Omega$,
one has
\begin{equation}\label{levi1} L(\varphi)(a)=\left|
\begin{matrix}
0 & \varphi_z & \varphi_w \\
\varphi_{\bar z} & \varphi_{z \bar z} & \varphi_{w \bar z} \\
\varphi_{\bar w} & \varphi_{z \bar w} & \varphi_{w \bar w} 
\end{matrix}
\right| \geq 0, \quad a \in \del \Omega,
\end{equation}
where $(z,w)$ are the variables of $\C^2$. For general $n \geq 2$,
J. Krzoska (1933) formulated it as with the same boundary
regularity, the hermitian matrix
\begin{equation}\label{levi2}
 \left(\frac{\del^2 \varphi}{\del z_j \del \bar z_k} (a)
\right)_{1 \leq j,k \leq n} \quad (a \in \del \Omega)
\end{equation}
is {\em positive semi-definite} on the homomorphic tangent vector space
\[
 \left\{(v_1, \ldots, v_n) \in \C^n: \sum_{j=1}^n v_j
 \frac{\del\varphi}{\del z_j} (a)=0 \right\}.
\]
If $n=2$, this is reduced to \eqref{levi1}.
Then it is natural to ask the converse.

{\bf Levi Problem:}  {\em If $\del \Omega$ satisfies \eqref{levi2}, is
$\Omega$ a domain of holomorphy?}

The property characterized by \eqref{levi1} or \eqref{levi2} is called
a {\em pseudoconvexity} of $\Omega$ or $\del \Omega$, which is
a biholomorphically invariant property
 in a neighborhood of any point $a \in \del \Omega$.

There is an inconvenience in the above characterization by $\varphi$;
that is, even if $\varphi_1, \varphi_2$ satisfies \eqref{levi1} or
\eqref{levi2}, $c_1 \varphi_1+c_2 \varphi_2$ with positive
constants $c_1, c_2$, does not satisfy the similar condition.
This was the reason why K. Oka introduced a {\em pseudoconvex} function
$\psi$ in $\Omega$ such that $\psi$ is upper semi-continuous and
the restriction of $\psi$ to the intersection of any complex
affine line and $\Omega$ is subharmonic
 (Oka VI, 1942).\footnote{\,In similar time,
 P. Lelong defined the same notion as
{\em plurisubharmonic} functions from potential theoretic viewpoint.}
Pseudoconvex functions play the similar role to that of
 $\varphi$ in  \eqref{levi1} or \eqref{levi2} and still satisfies that $c_1\psi_1+c_2 \psi_2$
is pseudoconvex for pseudoconvex functions $\psi_j$ and $c_j>0$
($j=1,2$).
If $\psi: \Omega \to \R$ is of $C^2$-class, $\psi$ is pseudoconvex
if and only if the hermitian matrix
$\left(\frac{\del^2 \psi}{\del z_j \del \bar z_k} (a)
\right)_{1 \leq j,k \leq n}$ ($a \in \del \Omega$) is positive
semi-definite.

In the unpublished papers 1943, K. Oka did not assume the
 boundary regularity of $\Omega$,
but defined the {\em pseudoconvexity} of $\Omega$ (or $\del \Omega$)
 as follows:
For every point $a \in \del \Omega$ there is a
neighborhood $U$ of $a$ in $\C^n$ such that
if $\phi: \Omega_{\mathrm{H}}(a;\delta, \gamma) \to U \cap \Omega$
is a biholomorphic map from
a Hartogs domain $\Omega_{\mathrm{H}}(a;\delta, \gamma)$
into $U \cap \Omega$, 
then $\phi$ is analytically continued to
$\tilde\phi: \PD(a; \gamma) \to U \cap \Omega$.
It is trivial that a domain of holomorphy satisfies this
pseudoconvexity, and K. Oka proved the converse:
This is why he called the problem {\bf Hartogs' Inverse Problem}.
The solution naturally implies that of the Levi Problem.

(ii)  Let $\Omega=\bigcup_{\alpha \in \Gamma} U_\alpha$ be an open
covering. Let $f_\alpha ~(\alpha \in \Gamma)$ be a meromorphic
function in $U_\alpha$ such that $f_\alpha - f_\beta$ is
holomorphic in $U_\alpha \cap U_\beta$ as far as
$U_\alpha \cap U_\beta \not= \emptyset$. The  pair
 $(\{U_\alpha\}, \{f_\alpha\})$ is called a Cousin-I data on $\Omega$.

{\bf Cousin I Problem}\footnote{\label{cousin}  This problem was dealt with
by P. Cousin \cite{co} and affirmatively solved when
 the domain is a {\em cylinder (domain)}\label{cyl},
which is by definition an $n$-product of the coordinate plane domains
of $\C^n$: Cousin II Problem below
was also solved affirmatively there
 when the domain is a cylinder $\prod_j D_j$ with
simply connected plane domains $D_j ~(\subset \C)$ except for one $D_j$.
He used the so-called {\em Cousin integral} (see  p.~\pageref{coin})
They were solved affirmatively in general by
 K. Oka \cite{oka}, I---III for univalent domains,
which in Cousin II Problem yielded the {\em Oka Principle}.}{\bf:}
 {\em If $\Omega$ is a domain of holomorphy,
then for a Cousin-I data  $(\{(U_\alpha\}, \{f_\alpha\})$ on $\Omega$,
find a meromorphic function $F$ in $\Omega$, called a solution
of the Cousin-I data,
such that $F - f_\alpha$ is holomorphic in every $U_\alpha$.}

In the case of $n=1$, Mittag-Leffler's Theorem gives
 an affirmative answer to the problem.

Similarly,  we assume that $f_\alpha$ are meromorphic functions,
not identically zero, and that $f_\alpha / f_\beta$ is a nowhere
vanishing holomorphic function in every
 $U_\alpha \cap U_\beta (\not=\emptyset)$.
 Then $(\{U_\alpha\}, \{f_\alpha\})$ is called a Cousin-II data
on $\Omega$.

{\bf Cousin II Problem:} {\em If $\Omega$ is a domain of holomorphy,
then for a Cousin-II data  $(\{U_\alpha\}, \{f_\alpha\})$ on $\Omega$,
find a meromorphic function $F$ in $\Omega$, called a solution
of the Cousin-II data,
such that $F/ f_\alpha$ is nowhere zero holomorphic
 in every $U_\alpha$.}

In the case of $n=1$, this is answered affirmatively
by Weierstrass' Theorem.

(iii)  Let $K \Subset \Omega$ be a compact subset and let $f$ be
a holomorphic function in a neighborhood of $K$.

{\bf Problem of expansion (Approximation Problem):} {\em
Assume that $\Omega$ is a domain of holomorphy.
Find a condition for $K$ such that for every such $f$ there is a series
 $\sum_{\nu=1}^\infty f_\nu$ with
holomorphic functions $f_\nu$ in $\Omega$ such that restricted on $K$,
\[
 f=\sum_{\nu=1}^\infty f_\nu,
\]
where the convergence is uniform on $K$.} 

In the case of $n=1$ we have Runge's Theorem.
In the problems of (ii) and (iii) above, the assumption for
$\Omega$ being a domain of holomorphy is necessary by examples
(cf., e.g., \cite{nog19b} \S1.2.4, \S3.7).

\section{Unpublished Papers VII---XI 1943}\label{7-11}
We first list the titles translated from Japanese
 and the numbers of pages of the papers with dates.
\begin{enumerate}
\item
On Analytic Functions of Several Variables VII ---
Two auxiliary problems on the congruence of holomorphic functions,
pp.~28 (4 Sep.\ 1943).
\item
On Analytic Functions of Several Variables VIII ---
The First Fundamental Lemma on finite domains without ramification
     points, pp.~11 (5 Sep.\ 1943).
\item
On Analytic Functions of Several Variables IX ---
Pseudoconvex functions, pp.~29 (24 Oct.\ 1943).
\item
On Analytic Functions of Several Variables X ---
The Second Fundamental Lemma, pp.~11 (12 Nov.\ 1943).
\item
On Analytic Functions of Several Variables XI ---
Pseudoconvex domains and  finite domains of holomorphy:
 Some theorems on finite domains of holomorphy, pp.~29 (12 Dec.\ 1943).
\end{enumerate}

K. Oka cited these papers in two places of the published papers
with mentioning a further problem of ramified Riemann domains,
which we quote.

(1) Introduction of  \cite{oka} Oka VIII (1951, p.\,204) begins with:
\vspace{-7pt}
\begin{quote}
  Les probl\`emes principaux depuis le
M\'emoire I sont : probl\`emes de Cousin, probl\`eme de d\'eveloppement 
et probl\`eme des convexit\'es\footnote{\,Ces probl\`emes sont fond\'es
sur H. Behnke et P. Thullen, Theorie der Funktionen mehrerer
Komplexer Ver\"anderlichen, 1934. Nous allons les expliquer en formes
pr\'ecises. Soient $\mathfrak{D},\mathfrak{D}_0$ deux domaines connexes
ou non sur l'espace de $n$ variables complexes tels que
$\mathfrak{D}_0\subseteq \mathfrak{D}$ (c'est-\`a-dire que
$\mathfrak{D}_0$ soit un $\ll$Teilbereich$\gg$ de $\mathfrak{D}$); nous
appellerons que $\mathfrak{D}_0$ est holomorphe-convexe par rapport
\`a $\mathfrak{D}$, s'il existe une fonction holomorphe dans
$\mathfrak{D}$ ayant des \'el\'ements de Taylor diff\'erents aux points
diff\'erents de $\mathfrak{D}_0$ et encore si, pour tout domaine connexe
ou non $\Delta_0$ tel que $\Delta_0 \Subset \mathfrak{D}_0$
(c'est-\`a-dire que $\Delta_0 \subset \mathfrak{D}_0$ et
$\Delta_0 \ll \mathfrak{D}_0$), on peut trouver un domaine connexe ou
non $\Delta$ tel que $\Delta_0\subseteq\Delta\Subset
\mathfrak{D}_0$ de fa\c con qu'\`a tout point $P$ de
$\mathfrak{D}_0-\Delta$, il corresponde une fonction $f$ holomorphe
dans $\mathfrak{D}$ telle que $|f(P_0)|>\max|f(\Delta_0)|$. 
Sp\'ecialement, si $\mathfrak{D}_0$ est ainsi par rapport \`a
lui-m\^eme, nous l'appelons avec H. Behnke d'\^etre holomorphe--convexe
(regul\"ar--konvex). Les probl\`emes sont alors : Probl\`emes de
Cousin. Trouver une fonction m\'eromorphe (ou holomorphe) admettant
les p\^oles (ou les z\'eros satisfaisant \`a une certaine condition)
donn\'es dans un domaine holomorphe--convexe. Probl\`eme de
d\'eveloppement. Soit $\mathfrak{D}_0$ un domaine (connexe ou non)
holomorphe--convexe par rapport \`a $\mathfrak{D}$; trouver, pour
toute fonction holomorphe $f$ une s\'erie de fonctions holomorphes 
dans $\mathfrak{D}$, convergente uniform\'ement vers $f$ dans tout
domaine connexe ou non $\Delta_0$ tel que $\Delta_0\Subset
\mathfrak{D}_0$. Probl\`eme des convexit\'es. Tout domaine
pseudoconvexe est-il holomorphe--convexe ? Pour les domaines
univalents, on peut remplacer $\ll$holomorphe-convexe$\gg$ par
$\ll$domaine d'holomorphie$\gg$, gr\^ace au th\'eor\`eme de H.\ Cartan
et P.\ Thullen.}.   
Dans les M\'emoires I--VI\footnote{\,Les M\'emoires pr\'ec\'edents sont :
I--Domaines convexes par rapport aux fonctions rationnelles, 1936; 
II--Domaines d'holomorphie, 1937; III--Deuxi\`eme probl\`eme de
Cousin, 1939 (Journal of Science of the Hiroshima University); 
IV--Domaines d'holomorphie et domaines rationnellement convexes, 1941; 
V--L'int\'egrale de Cauchy, 1941 (Japanese Journal of Mathematics); 
VI--Domaines pseudoconvexes, 1942 (Tohoku Mathematical Journal); VII--Sur
quelques notions arithm\'etiques, 1950 (Bulletin de la Soci\'et\'e
Math\'ematique de France)}, nous avons vu, disant un mot, que ces
probl\`emes sont r\'esolubles affirmativement pour les domaines
univalents finis\footnote{\,Pr\'ecis\'ement dit, pour le deuxi\`eme
probl\`eme de Cousin, nous avons montrer une condition n\'ecessaire et
suffisante pour les z\'eros; et pour le probl\`eme des convexit\'es,
nous l'avons expliqu\'e pour les deux variables complexes, pour diminuer 
la r\'ep\'etition ult\'erieure in\'evitable.}.  
Et l'auteur a encore constat\'e quoique sans l'exposer, que ces
r\'esultats restent subsister au moins jusqu'aux domaines finis sans
point critiques\footnote{\,L'auteur l'a \'ecrit aux d\'etails en
japonais \`a Prof. T. Takagi en 1943.}.

Il s'agit donc: ou bien d'introduire l'infini convenable, ou bien de
permettre des points critiques; or, on retrouvera que l'on ne sais
presque rien sur les domaines int\'erieurement ramifi\'es; .....
\end{quote}

(2) Introduction 2 of \cite{oka} Oka IX (1953, p.\,98) begins with:
\vspace{-7pt}
\begin{quote}
 Dans le pr\'esent M\'emoire, nous traiterons les
probl\`emes indiqu\'es plus haut, ainsi que les probl\`emes
arithm\'etiques introduits au M\'emoire VII, pour les domaines
pseudoconvexes finis sans point critique int\'erieur; dont la partie
essentielle n'est pas diff\'erente de ce que nous avons expos\'e en
japonais en 1943\footnote{\,Voir la Note \`a l'Introduction de M\'emoire 
VIII. Dans ce manuscrit-ci on trouve d\'ej\`a les probl\`emes
$(\mathrm{C_1})\ (\mathrm{C_2})$ (expricitement) et $(\mathrm{E})$
(implicitement).}.     

On verra dans le M\'emoire suivant que quand on admet les points
critiques int\'erieurs, on rencontre \`a un probl\`eme qui
m'appara\^\i t extr\^emement difficile (voir No.\ 23). C'est pour
pr\'eparer des m\'ethodes et pour \'eclaircir la figure de la
difficult\'e, que nous avons d\'ecid\'e \`a publier le pr\'esent
M\'emoire, s\'epar\'ement\footnote{\,cit\'e plus haut.}.  
\end{quote}

For convenience we recall their English translations
by R. Narasimhan from \cite{oka2}:

\begin{itemize}
\item[(1)] The principal problems we have dealt with since Memoir~I are
the following: Cousin problems, the problem of expansions and
 the problem of (different types of)
 convexity\footnote{\,These problems are based on H. Behnke
 and P. Thullen, Theorie der Funktionen mehrerer
komplexer Ver\"anderlichen, 1934. Let us explain them in precise form.
 Let $\mathfrak{D}, \mathfrak{D}_0$ be two domains over the space of
 $n$ complex variables connected or not such that
$\mathfrak{D}_0 \subseteqq \mathfrak{D}$
(i.e.\ such that $\mathfrak{D}_0$ is a ``Teilbereich'' of $\mathfrak{D}$). 
We shall say that $\mathfrak{D}_0$ is holomorph-convex with respect to $\mathfrak{D}$
if $\mathfrak{D}_0 \subseteqq H$, $H$ being the ``Regularit\"atsh\"ulle''
of $\mathfrak{D}_0$, and if, in addition, for every domain $\Delta_0$,
connected or not, such that $\Delta_0  \Subset \mathfrak{D}_0$
(that is, $\Delta_0 \subset \mathfrak{D}_0$ and $\Delta_0 \ll \mathfrak{D}_0$),
we can find a domain $\Delta$, connected or not such that
$\Delta_0 \subset \Delta \Subset \mathfrak{D}_0$ and such that,
to every point $P_0$ of $\mathfrak{D}_0-\Delta$, there
corresponds a function $f$ holomorphic on $\mathfrak{D}$ with
$f(P_0) > \max | f(\Delta_0) |$. In particular, if $\mathfrak{D}_\alpha$ has
this property with respect to itself, we call it, 
 with H. Behnke, holomorph-convex (regul\"arkonvex).
 The problems are then the following: Cousin problems.
 Find a meromorphic (or holomorphic) function having given poles
 (or given zeros satisfying a certain additional condition).
 Problem of expansions. Let $\mathfrak{D}_0$
 be a domain (connected or not) holomorph-convex
with respect to $\mathfrak{D}$; for any function f holomorphic on $\mathfrak{D}_0$,
 find a series of holomorphic functions on $\mathfrak{D}$ 
which converges uniformly to $f$ on any domain $\Delta_0$,
connected or not, such that $\Delta_0 \Subset \mathfrak{D_0}$.
Problem of convexity. Is every pseudoconvex domain holomorph-convex?
 For univalent domains, one can replace
``holomorph-convex'' by ``domain of holomorphy'' because of
the theorem of H. Cartan and P. Thullen.}
 In Memoirs I---VI\footnote{\,The preceding Memoirs are:
 I. Rationally convex domains, 1936; II. Domains of holomorphy, 1937;
III. The second Cousin problem, 1939 
(Journal of Science of Hiroshima University);
 lV. Domains of holomorphy and rationally convex domains, 1941; 
V. The Cauchy integral, 1941 (Japanese Journal of Mathematics);
VI. Pseudoconvex domains, 1942 (Toh\^oku Mathematical Journal);
VII. On some arithmetical concepts, 1950 (Bulletin de la
 Soci\'et\'e Mathematique de France)}
we have seen, to put it in one word, that these problems can 
be solved affirmatively for univalent domains
without points at infinity\footnote{\,More precisely,
 we obtained a necessary and sufficient condition for the second Cousin
problem; and the problem of convexity was only explained
 for two complex variables in order to reduce the ultimate
 repetition which is inevitable.}.
 Furthermore, the author has verified, albeit without publishing this,
 that these results remain valid at least as far as
domains without points at infinity and without
 interior ramification points\footnote{\,The author
 has written this out in detail in Japanese and sent it
 to Prof. T. Takagi in 1943.}.

We must therefore either introduce suitable points at infinity or allow
points of ramification. Now, one will find that almost nothing is known about
domains with interior ramification. .....

\item[(2)]
 In the present memoir, we shall deal with the problems indicated above,
as well as the arithmetical problems introduced in Memoir VII,
for pseudoconvex domains without interior ramification and
 without points at infinity; the essential part of this memoir
 is not very different from what we have expounded
 in Japanese in 1943\footnote{\,See the note in the introduction 
to Memoir VIII. In that manuscript, one finds already
problems ($C_1$),  ($C_2$) (explicitly), and problem (E) (implicitly).
}.

We shall see in the memoir following this one that when one permits
interior points of ramification, one meets a problem which seems 
to me to be extremely difficult (see also No.~23 below).
 It is to prepare the methods and to illuminate the nature of
 this difficulty that we have decided to publish the
present memoir separately\footnote{\,Cite the above.}.
\end{itemize}

According to T. Nishino (\cite{okap} Vol.~1, Afterword),
the original manuscripts of this series sent to T. Takagi in 1943 were
lost, but fortunately, the complete set of their draft-manuscripts
had been kept in Oka's home library and was found posthumously.

It is really surprising for me to learn that the way of arguments
in Oka IX (published, 1953) is very similar to the one in the
series of papers 1943, ten years prior, and that the part of the arguments
to prove so-called {\em Oka's Heftungslemma}\footnote{\,Roughly speaking,
the union of two adjacent holomorphically convex domains
with pseudoconvex boundary is holomorphically
convex (cf., e.g., \cite{an})},
 an essential step in the proof of the Levi (Hartogs' Inverse) Problem,
 is almost a copy of the corresponding part in unpublished Paper XI 1943.

\smallskip
For the English translation of Paper XI, I describe in below
some supplements and recall briefly the main results
that had been obtained in VII---X and used in XI.

H. Cartan once has written (\cite{oka2}, p.\,XII):\label{carw}
\vspace{-10pt}
\begin{quote}\sf
.............\\
\quad Mais il faut avouer que les aspects techniques de ses d\'emonstrations
et le mode de pr\'esentation de ses r\'esultats rendent difficile la 
t\^ache du lecteur, et que ce n'est qu'au prix d'un r\'eel effort que
l'on parvient \`a saisir la port\'ee de ses r\'esultats,
qui est consid\'erable. 
C'est pourquoi il est peut-\^etre encore utile aujourd'hui, en hommage
 au grand cr\'eateur que fut Kiyoshi O{\sc ka}, de pr\'esenter
l'ensemble de son {\oe}uvre.\\
\quad .................
\end{quote}
In English (by Noguchi),
\vspace{-10pt}
\begin{quote}\sf
.............\\
\quad But 
we must admit that the technical aspects of his proofs
and the mode of presentation of his results make it difficult
to read, and that it is possible only at the cost of a real effort to
grasp the scope of its results, which is considerable.
This is why it is perhaps still useful today,
for the homage of the great creator that was Kiyoshi O{\sc ka},
 to present the collection of his work.\\
\quad .................
\end{quote}

The present series is no exception.
The aim of the series is two folded:
\begin{enumerate}
\item
With an intention to deal with the problem for {\em ramified}
 Riemann domains, the
     conditions and the statements of lemmata, propositions etc.\
 are made as general as possible.
\item
In the same time, they must be satisfied and proved completely
for {\em unramified} Riemann domains as a special case.
\end{enumerate}
This approach which contains in a sense a self-confliction between
``general'' versus ``special'' seems to increase an involvedness
of the presentations of the papers, but forms a motivation
 to invent ``Coherence'' or 
{\em ``Id\'eaux de domaines ind\'etermin\'es}''
in terms of Oka (see \S\ref{43-}), which is referred, e.g.,
as:
\begin{quote}{\sf
Of greatest importance in Complex Analysis is the concept of
a coherent analytic sheaf} (Grauert--Remmert \cite{grcas}).
\end{quote}

The last Paper XI contains the most important final conclusion proving
that {\em every pseudoconvex unramified Riemann domain over $\C^n$
is Stein (in terms of the present days)}.
In a year before, 1942, Oka published Oka VI (\cite{oka}),
proving the result in the case of univalent domains
 of $\C^2$.
In Oka VI (\cite{oka}),
 he used Weil's integral formula, which in $n$-dimensional
case takes a rather involved form already in univalent domains.
To deal with possibly infinitely sheeted unramified Riemann domains
with his intension even to deal with ramified case, he wanted
to  avoid the use of Weil's integral formula, but to use simpler
Cauchy's integral formula combined with ``{\em J\^oku-Ik\^o\,}''
(lifting principle) which was prepared
as {\em The First Fundamental Lemma} at the end of
 Paper VIII of the present series.
 The method of J\^oku-Ik\^o was invented
in his first two papers Oka \cite{oka} I and II.

For the proof of ``Heftungslemma'',
he uses an integral equation of the Fredholm type
similarly to Oka VI (published, \cite{oka}); in Oka IX
 (published, \cite{oka}) the integral equation is implicit.

Reading the series of unpublished papers VII---XI 1943,
 we observe not only the solution of
 the Levi (Hartogs' Inverse) Problem for unramified
Riemann domains over $\C^n ~(n \geq 2)$,
but also  the dawn of the
then unknown notion of {\em ``Id\'eaux de domaines ind\'etermin\'es}''
or ``Coherence''.

\begin{rmk}\rm It is a nature of Oka's wording such as
{\em Id\'eaux de domaines ind\'etermin\'es} to represent
``{\em a way of thinking}'' rather than the formed object,
 similarly to the case of
``{\em J\^oku-Ik\^o}'' (see Footnote \ref{jkik} at p.\,\pageref{jkik}).
\end{rmk}

\section{The XI-th paper}
\subsection{Some practical notes}
This series of the present Papers VII---XI in 1943
 were written as a continuation of the
published papers Oka I---VI (\cite{oka}).
In Part II we shall present a Japanese translation
of the last Paper XI, in which 
at some important places, footnotes are put to remind the numbering
as ``Note by the translator''.
As a consequence, the numbering of the footnotes are different
to the original.

As Oka writes ``Report VI'', then it means the published paper with
the same number in \cite{oka}. On the other hand, Report VII to X
(e.g., Report IX) is the article of the present series
(not the published Oka IX in \cite{oka}).

As Oka writes ``a {\em finite} domain'', it means a multi-sheeted
domain spread over $\C^n$, not containing an infinite point,
say, in a compactification such as complex projective $n$-space.

\subsection{The XI-th paper}
This is the last one of the series from VII-th, in which
Oka settled affirmatively
the Levi (Hartogs' Inverse) Problem for general dimensional
unramified Riemann domains over $\C^n$, ten years
before Oka \cite{oka} IX was published in 1953:
There was then no notion of
``{\em Coherence}'' or ``{\em Id\'eaux de domaines ind\'etermin\'es}''
termed by Oka.
It is rather surprising to know that the Problem had been solved
just after Oka VI 1942 (in the case of
$2$-dimensional univalent domains) by a different method,
if one observes the state of advances at that time as discussed
in \S\ref{3bp}.

Because of the importance, I chose the last one for the
translation into English.

In this paper K. Oka begins with proving the {\em Cousin I/II Problems}
as well as the {\em Problem of expansions (Approximation Problem)
for unramified Riemann domains over $\C^n ~(n \geq 2)$}
by a different method than those in Oka \cite{oka} I---III,
using a new J\^oku-Ik\^o prepared in Papers VII--VIII
of the present series.

Let us quote the most important main result from Paper XI \S10 (Part II):
\begin{quote}\sf
{\bf Theorem I.}  {\em
A finite pseudoconvex domain with no interior ramification point is
 a domain of holomorphy.}
\end{quote}

\begin{rmk}\rm
\begin{enumerate}
\item
In the published Oka I--VI the domains are assumed to be univalent.
Oka first dealt with  unramified multivalent domains over $\C^n$ systematically
in the present series of VII--XI.
\item
In the proof of Oka's Theorem I above he actually proves that such a pseudoconvex
     domain is holomorphically convex and satisfies the separation
property by holomorphic functions (see the footnote of Theorem I, XI \S10).
It is noted that unramified holomorphically convex domains (multivalent
     in general) are domains of holomorphy; the converse holds,
provided that the domains are finitely sheeted (due to Cartan--Thullen
     \cite{ct}). Cartan--Thullen \cite{ct} claimed the converse in
     general, but there was an oversight in the case of infinitely many sheeted domains.
The oversight was fulfilled by the proof of Oka's Theorem I above (cf.\ XI \S11)
 as a series of implications:
``domain of holomorphy'' $\Rightarrow$ ``pseudoconvex domain''
$\Rightarrow$ ``holomorphically convex domain''.
Thus the three classes of unramified domains over $\C^n$ are equivalent.
\end{enumerate}
\end{rmk}

\section{The VII---X-th Papers}

To begin with,
 it will be interesting and worthy to recognize Oka's own observation
of the state of researches at the time to start writing
 the present series of papers 1943
 by recalling the first paragraph of the VII-th:
\begin{quote}\sf
The problems discussed at the beginning of the first report\footnote{\,
(Note by the present author)
This is Oka \cite{oka} I 1936; the same in the sequel.}
 were solved more or less generally in the series of reports up to VI\footnote{\,
(The original footnote) I began the present research with the back
 ground of the following monograph.\\
H. Behnke--P. Thullen: Theorie der Funktionen mehrerer komplexer
 Ver\"anderlichen, 1934 (Egebnisse der Mathematik und ihrer
 Grenzgebiete).

Reports prior to this report: I, 1936; II, 1937; III, 1939,
(Journal of Science of the Hiroshima University).
IV, 1940; V, 1940, (Japanese Journal of Mathematics).
 VI, 1942, (The T\^ohoku Mathematical Journal).
 }.
But, since these were a sort of depth sounding in a sense,
 we avoided domains such as not finite or non univalent,
 and considered some of them only in
the case of two variables. While we may think of
 really various kinds of problems
on analytic functions of several variables, it is, for a moment,
our main aim of the research  to get rid of these restrictions one
 by one. The present paper is devoted to the preparation for it.
\end{quote}

In fact, Oka begun in the present series of papers
to deal with  general multi-valent domains over $\C^n$, systematically.
Here we would like to summerize briefly what were proved
in the VII---X-th papers before the XI-th paper.

The four papers were roughly classified into two groups,
VII+VIII and IX+X.

\subsection{VII+VIII}\label{jokuiko}
These two papers were devoted to the study of ideal theoretic
properties of holomorphic functions. The study of this part
led to the works of ``{\em Id\'eaux de domaines ind\'etermin\'es}''
or ``{\em Coherence}'' (Oka VII, VIII, published \cite{oka}).
Therefore, in Oka IX (published, \cite{oka}) the contents of
this part were replaced by the more general results of
 Oka VII, VIII (published, \cite{oka}).

In VII he considered a domain $\mathfrak{D}$ in the space of $n$
complex variables $x_1, \ldots, x_n$. Let $\O(\mathfrak{D})$ denote
the ring of all holomorphic functions in $\mathfrak{D}$.
Let $(F)=(F_1, F_2, \ldots, F_p)$ be a system of holomorphic functions
in $\mathfrak{D}$.
For $f(x), \varphi(x) \in \O(\mathfrak{D})$ we write
\[
 f \equiv \varphi \quad (\mathrm{mod.}~ F_1, F_2, \ldots, F_p),
\]
and say that $f$ and $\varphi$ are {\em congruent} with respect to
the function system $(F)$ in $\mathfrak{D}$,
if there are functions $\alpha_j \in \O(\mathfrak{D})$ ($1 \leq j \leq p$)
satisfying
\[
 f-\varphi=\alpha_1 F_1 +\alpha_2 F_2 + \cdots + \alpha_p F_p.
\]

Let $P$ be a point of $\mathfrak{D}$. We define the notion of
{\em being congruent at $P$} if the above property hold in a
neighborhood of $P$. Then it is different to say that they are
congruent in $\mathfrak{D}$ and they are congruent at each point
of $\mathfrak{D}$. To emphasize this difference we also say
the former case to be {\em congruent globally in $\mathfrak{D}$}.

If $\bar{\mathfrak{D}}$ is a closed domain, we denote by
$\O(\bar{\mathfrak{D}})$ the set of all of holomorphic functions
in neighborhoods of $\O(\bar{\mathfrak{D}})$.

Then he formulate two problems:

{\bf Problem I.} {\em Let $\bar{\mathfrak{D}}$ be a bounded closed domain
in $(x)$ space.
 For a given  holomorphic function system
$(F)=(F_1, F_2, \ldots, F_p)$ with $F_j \in \O(\bar{\mathfrak{D}})$
and a given holomorphic function $\Phi(x) \in \O(\bar{\mathfrak{D}})$
such that  $\Phi(x)\equiv 0 ~ (\mathrm{mod.} ~ F)$
at every point $P \in \bar{\mathfrak{D}}$,
choose $A_j \in \O(\bar{\mathfrak{D}})$ so that}
\[
 \Phi(x)=A_1(x)F_1(x)+A_2(x)F_2(x)+ \cdots + A_p(x)F_p(x),
\quad x \in \bar{\mathfrak{D}}.
\]

{\bf Problem II.} {\em Let $(F)=(F_1, F_2, \ldots, F_p)$ be
a system of holomorphic functions defined in a neighborhood of
 $\bar{\mathfrak{D}}$. Suppose that for each point
 $P \in \bar{\mathfrak{D}}$ there are associated a
polydisk $(\gamma)$ with center $P$ and a holomorphic
function $\varphi(x)$ in $(\gamma)$ satisfying that
for two such pairs $((\gamma_j), \varphi_j), j=1,2$, with
$(\delta)=(\gamma_1) \cap (\gamma_2)\not=\emptyset$,
\[
 \varphi_1(x)\equiv \varphi_2(x) \quad (\mathrm{mod. }~
F_1, F_2, \ldots, F_p)
\]
at every point of $(\delta)$ (congruent condition).
Then, find a $\Phi(x) \in \O(\bar{\mathfrak{D}})$ such that
\[
 \Phi(x) \equiv \varphi(x) \quad (\mathrm{mod. }~ F)
\]
at every point $P \in \bar{\mathfrak{D}}$.}

\begin{rmk}\rm
Problem I is a sort of Syzygy type problem, and
Problem II is a Cousin-I Problem for the ideal generated
by $(F)=(F_1,F_2, \ldots, F_p)$.
\end{rmk}

In \S2 of Paper VII he defines the following property named

$(A)$:
 {\em Let $(F_1, F_2, \ldots, F_p)$ be a system of holomorphic
functions in a domain $\mathfrak{D}$ of $(x)$-space such that
$F_1 \not \equiv 0$. Let  $q \in \{2,3, \ldots, p\}$
and let  $P \in \mathfrak{D}$ be an arbitrary point.
If holomorphic functions $\alpha_j(x) ~ (j=1,2, \ldots, q)$ in
a neighborhood $U (\subset \mathfrak{D})$ of $P$ satisfy
\[
 \alpha_1(x)F_1(x)+ \alpha_2(x)F_2(x)+ \cdots + \alpha_q(x)F_q(x)=0,
~~ x \in U,
\]
then
\[
 \alpha_q(x) \equiv 0 \quad (\mathrm{mod.}~  F_1, F_2, \ldots, F_{q-1})
\quad \hbox{at }~ P.
\]
}

Most importantly, he shows 
the following for property (A):

\smallskip
{\bf Lemma 1.}  {\em Let $X$ be a domain in $(x)$-space, and let
$f_j(x) ~(j=1,2, \ldots, \nu)$ be holomorphic functions in $X$.
 Then the system
of holomorphic functions $F_j(x,y)=y_j -f_j(x)$ ($j=1,2, \ldots, \nu$)
satisfies property $(A)$.}

This is intended to apply for an Oka map
\[
 \psi(x)=(x, f_1(x), f_2(x), \ldots , f_\nu(x)) \in
\Omega \times \Delta(1)^\nu \subset \Delta(R)^n \times \Delta(1)^\nu,
\]
where $f_j(x) \in \O(X)$, $\Omega ~(\Subset X)$ is an analytic
polyhedron defined by
\[
x \in X, \quad |f_j(x)|<1, \quad j=1,2, \ldots, \nu,
\]
{$\Delta(R)$ is the disk of radius $R ~(>0)$ with center at the origin
in $\C$ and $R$ is chosen so that $\Omega \subset \Delta(R)^n$.
This is the essential part of Oka's {\em J\^oku-Ik\^o}:}

\medskip
{{\bf Remark (J\^oku-Ik\^o).} T. Nishino \cite{nis96} uses
``lifting principle'' for ``J\^oku-Ik\^o''.
It is a methodological principle termed by Oka such that}
\begin{enumerate}
\item
{one embedds a domain
into a higher dimensional domain of simple shape
(i.e., a polydisk) through the Oka map above;}
\item
{one extends a difficult problem on the original domain to the one on
the higher dimensional domain of simple shape;}
\item
{by making use of the simpleness of the higher dimensional domain,
one obtains a solution of the problem;}
\item
{then, one restricts the solution on the embedded original domain
to get a solution of the original problem.}
\end{enumerate}
{Things do not go so simply,
 but this is the principal method of K. Oka
all through his works.}

\medskip
Oka then affirmatively solves  Problems I and II under
 this property $(A)$ for $(F)$.

\medskip
\textbf{Theorem 1.} {\em
Let  $\bar{\mathfrak{D}}$ be a bounded
 closed cylinder domain
and let $(F)=(F_1, F_2, \ldots, ,F_p)$ be a system of
 holomorphic functions in a neighborhood of
$\bar{\mathfrak{D}}$ which  satisfies property $(A)$.
Then, Problem I for $(F)$ is solvable.} 

For a cylinder domain, see  Footnote~\ref{cousin}, p.~\pageref{cyl}.

\medskip
\textbf{Theorem 2.} {\em
Let  $\bar{\mathfrak{D}}$ and $(F)$ be the same as in Theorem 1 above.
Then, Problem II for $(F)$ is solvable.} 

In \S\S8--10 of Paper VII Oka deals with Problems I and II with estimates.

Finally, at the end of Paper VIII Oka obtained

\smallskip
{\bf Fundamental Lemma I. }  {\em Let $X$ be a univalent cylinder domain
in $(x)$-space and $\Sig \subset X$ be an analytic subset.
Let $V$ be a univalent open subset of $X$, containing $\Sig$.
Suppose that there are holomorphic functions
$f_1(x), f_2(x), \ldots, f_p(x) \in \O(V)$ such that
$\Sig=\{x \in V: f_j(x)=0, 1 \leq j \leq p\}$.
Let $X^0 \Subset X$ be a univalent bounded cylinder domain,
and set $\Sig_0=\Sig \cap X^0$.

Then, for every $\varphi(x) \in \O(V)$ with
$|\varphi(x)|<M$ in $V$, there is a holomorphic function
$\Phi(x) \in \O(X^0)$ such that at every point of $\Sig_0$
\[
 \Phi(x) \equiv \varphi(x) \quad (\mathrm{mod. } f_1, f_2, \ldots, f_p),
\]
and
\[
 |\Phi(x)| < KM, \quad x \in X^0,
\]
where $K$ is a positive constant independent from $\varphi(x)$.
}

He finishes Paper VIII with writing
\begin{quote}\sf
This theorem should be generalized soon later, but
so far as we are concerned with finite domains without ramification
 points,
this is sufficient for our study.
\end{quote}

\begin{rmk}\rm
By this comment we see that he had in mind a project to
deal with Levi (Hartogs' Inverse) Problem generalized
to domains with ramifications.
\end{rmk}

\subsection{IX+X}

In these two papers Oka defines and studies
 {\em pseudoconvex}
functions, equivalently plurisubharmonic functions as well
strongly pseudoconvex (plurisubharmonic) functions, and
investigates the boundary problem of pseudoconvex domains.
The contents of these IX and X correspond to and appear in
Oka IX (published, \cite{oka}), Chap.~2, {\S\S}B and C.

In these papers he deals with domains, finite and unramified over
$(x)$-space of $n$ complex variables
$x_1, x_2, \ldots, x_n$.
He begins with the notion of unramified domains over $(x)$-space.

Let $\mathfrak{D}$ be a domain over $(x)$-space and let
 $E \subset \mathfrak{D}$ be a subset.
If the infimum of the Euclidean distances from $P \in E$
to the (ideal) boundary of $\mathfrak{D}$ is not $0$,
one says that {\em $E$ is bounded with respect to $\mathfrak{D}$}.

He defines a {\em pseudoconvex domain} modeled after
F. Hartogs as follows:

{\em Definition.} {A domain $\mathfrak{D}$ over $(x)$-space
is said to satisfy {\em Continuity Theorem} if the following
condition is satisfied: Let $r=(r_j), \rho=(\rho_j)$ be $n$-tuples
of positive numbers with $\rho_j < r_j$, and consider a
polydisk $\PD(a; r)$,  $|x_j -a_j|< r_j$ with center $a=(a_j)$ and
a Hartogs domain:
\[
 \begin{array}{lllrl}
 \Omega_{\mathrm{H}}(a; r, \rho): && |x_j - a_j|<\rho_j, & |x_n-a_n|<r_n 
 & (j=1,2,\ldots,n-1), \\
 & {\text{or}} 
 & |x_j - a_j|<r_j , & \rho_n<|x_n-a_n|<r_n    
 & (j=1,2,\ldots,n-1).
\end{array}
\]
If $\phi: \Omega_{\mathrm{H}}(a; r, \rho)  \to \mathfrak{D}$
is a biholomorphic map, then $\phi$ necessarily extends biholomorphically
to $\tilde\phi: \PD(a; r) \to \mathfrak{D}$.
}

{\em Definition.} A domain $\mathfrak{D}$ over $(x)$-space
is said to be {\em pseudoconvex} if
 the following two conditions are satisfied:
\begin{enumerate}
\item
For each boundary point $M$ of $\mathfrak{D}$ there is a positive
number $\rho_0$ with polydisk $\PD$ 
 of radius $\rho_0$
and center $\underline{M}$ of the underlying point of $M$
such that the maximal subdomain $\mathfrak{D}_0$ of $\mathfrak{D}$
with the boundary point $M$
whose  underlying points are contained in $\PD$ 
satisfies Continuity Theorem.
($\mathfrak{D}$ satisfies locally Continuity Theorem.)
\item
Let $\PD_1 \subset \PD$ 
 be a polydisk with the same
 center\footnote{\,The radius of each variable may different.},
and let $\mathfrak{D}_1$ be the maximal subdomain with
the boundary point $M$ whose underlying points are
 contained in $\PD_1$. 
Let $(T)$ be a  one-to-one quasi-conformal\footnote{\,It is
 unclear very much what ``quasi-conformal'' amounts to,
but it is {\em holomorphic}.}  transform
from $\PD_1$ 
into $(x')$-space with the image denoted by
 $\Delta'_1$, 
 and $\mathfrak{D}'_1=T(\mathfrak{D}_1)$.
Then, $\mathfrak{D}'_1$ satisfies always Continuity Theorem.
(The property (i) is not lost by quasi-conformal transforms.)
\end{enumerate}

\begin{rmk}\rm
From the definition  above one 
sees why he called the problem as {\em Hartogs' Inverse
Problem}.
\end{rmk}

Then he defines a {\em pseudoconvex function} or a
{\em plurisubharmonic function} valued in $[-\infty, \infty)$
 so that it is
upper-semicontinuous and its restriction to every complex
line segment is subharmonic.

After Hartogs' holomorphic radius he defines
 Hartogs' radii $R_j(P)$ ($j=1,2, \ldots, n$)
 at $P \in \mathfrak{D}$ by 
\[
 R_j(P)= \sup\{r_j: \PD(P; (r_1, \ldots, r_n)) \subset \mathfrak{D}, ~
r_j>0, ~ 1 \leq j \leq n \},
\]
where $\PD(P; (r_1, \ldots, r_n)) :=\{(z_j) \in \C^n: |z_j- p_j|<r_j,
1 \leq j \leq n\}$ with $P=(p_j)$.
He proves:

\smallskip
{\bf Theorem 1.} {\em If $\mathfrak{D}$ is pseudoconvex,
then $-\log R_j(P)$ is pseudoconvex in $\mathfrak{D}$.
(Here the logarithm stands for the real branch.)
}

Similarly, let $d(P)$ ($P \in \mathfrak{D}$) denote the supremum
of radii $r>0$ such that a  ball with center $P$ and radius $r$
is contained in $\mathfrak{D}$, and $d(P)$ is called the
{\em Euclidean boundary distance}.
He then proves:

\smallskip
{\bf Theorem 3.} {\em If $\mathfrak{D}$ is pseudoconvex,
then $-\log d(P)$ is a pseudoconvex function in $\mathfrak{D}$.
}

Then he consider a  pseudoconvex function
$\varphi(x)$  of $C^2$-class in general, confirming the semi-positivity
of the Hermitian form
\[
{W}(\varphi;  (v_j), (w_k))(P) = \sum_{j,k} \frac{\del^2 \varphi}
{\del x_j \del \bar{x}_k}(P) v_j \bar{w}_k, \quad
(v_j), (w_j) \in \C^n.
\]
This form $W(\phi; \cdot, \cdot)$, which was written so
in the paper and is nowadays
 called the {\em Levi form}, is due to Oka \cite{oka} VI.
Then he proves in IX:

\smallskip
{\bf Theorem 5.} {\em
If ${W}(\varphi; (v_j), (w_k))(P)$ is strictly positive definite
at $P=P_0$, then one can find a holomorphic
 polynomial function $f(x_1, x_2, \ldots, x_n)$ of degree $2$
such that $f(P_0)=0$ and in a neighborhood of $P_0$,
the analytic hypersurface $\{f=0\}$ lies in the
part $\{\varphi >0\}$ except for $P_0$.}

\begin{rmk}\rm
In one variable, the situation is much simpler:
If $\mathfrak{D}$ is a domain in $\C$ and $P_0 \in \del \mathfrak{D}$, then
$f(z)=z-P_0$. It is the purpose to construct
a meromorphic function on $\mathfrak{D}$ such that its poles
are only $\frac1{f(z)}$ near $P_0$.
 When $n \geq 2$,
Oka formulated the positivity of $W(\varphi; \cdot, \cdot)$
to have $f(z)$. Later, he solves the Cousin I Problem
on $\bar{\mathfrak{D}}$ with poles only $\frac1{f(z)}$ near $P_0$,
 and then concludes that
$\mathfrak{D}$ is holomorphically convex.
\end{rmk}

Oka took a smoothing of a pseudoconvex function $\varphi(x)$
by the volume integration average, and repeat it to have
a $C^2$-differentiable pseudoconvex function; nowadays it is more common
to take a convolution integration, but the role is the same.

Finally at the end of Paper X, he obtained

\smallskip
{\bf Fundamental Lemma II.} {\em
Let $\mathfrak{D}$ be a pseudoconvex domain over $(x)$-space
without ramification point. Then there is a continuous
pseudoconvex function $\varphi_0(P)$ in $\mathfrak{D}$
 satisfying the following two conditions:
\begin{enumerate}
\item
If $\mathfrak{D}_c:=\{P \in \mathfrak{D}: \varphi_0(P) <c\}$
for every real number $c$, then $\mathfrak{D}_c \Subset \mathfrak{D}$.
\item
There are exceptional points of $\mathfrak{D}$ with no
accumulation point inside $\mathfrak{D}$ and
for any other point $P_0 \in \mathfrak{D}$ than them,
 one can find an analytic hypersurface $\Sig$ passing
$P_0$ in a neighborhood of $P_0$ such that
$\varphi_0(P)> \varphi_0(P_0)$ for
$P \in \Sig \setminus \{P_0\}$.
\end{enumerate}
}

\section{After Paper XI, and Problem left}\label{43-}

The series of Papers VII--XI in 1943 was not translated into
French for publication, but continued to
Report XII dated 26 May 1944 (\cite{okap}), titled
\begin{itemize}
\item
On Analytic Functions of Several Variables ~XII ---
Representation of analytic sets, pp.~22.
\end{itemize}

In this manuscript, he first used {\em Weierstrass' Preparation Theorem}
to study  local properties of analytic sets.
As known well, Weierstrass' Preparation Theorem plays a
crucial role in the proofs of Oka's Coherence Theorems.
In this sense, 
the turn of years {\em 1943/'44} was indeed the {\em ``watershed''} in
the study of analytic function theory of several variables. 

The research was continued more to the following and further (cf.\ \cite{okap}):
\begin{itemize}
\item
XIII~ On the condition of Weierstrass' Preparation Theorem
(\cite{okap}).
\end{itemize}
The precise date of this manuscript is unclear, but probably
 around 16 November 1945 due to T. Nishino's comment
 (\cite{okap} Vol.\ 2, Afterword).

\begin{rmk}\label{44/45}\rm
These manuscripts are not completed ones, while the formers
 (VII---XI) are.
But it is still interesting to read the introductions of the above two
 manuscripts.
\begin{enumerate}
\item
In XII he begins with ``{\sf We wish to {\em extend more} the results obtained
in the former reports}''. 
\item
In XIII, he writes at the beginning: ``{\sf To consider the series of problems mentioned
at the beginning of this research, we have imposed some conditions
to the domains. From Reports I to VI\footnote{\,
(Note by the present author.) These numbers
 correspond to the published papers \cite{oka}.}  the domains were assumed to be
finite and univalent, and from VII to XI\footnote{\,
(Note by the present author.) These are the unpublished papers 1943.}
 we excluded the infinity points and interior ramification points from the domains.
As we are going to take out these conditions,
 it would be better to advance step by step.
Putting the problem of the infinity points aside for a moment,
 we firstly would like to investigate
what will happen if the {\em ramification points} are allowed
to the interior of the domains.}''\\
\quad It should be noticed that this paragraph was rephrased in the
published \cite{oka} VIII, 1951.
\end{enumerate}
\end{rmk}

\begin{rmk}\label{infram}\rm
The problem of the infinity points was affirmatively solved by
R. Fujita \cite{fuj} and A. Takeuchi \cite{tak} for unramified
domains over complex projective spaces with at least one ideal boundary
 point; there are more extensions in the unramified case, but we stop to
go further in this direction, which is away from the problem
of the ramified case.
A counter-example of a ramified domain over complex projective space
was given by H. Grauert \cite{nar}.
Therefore the problem of ramification points remained open
then for domains over $\C^n$: Later, in 1978 J.E. Forn{\ae}ss \cite{forn}
gave a counter-example of a two-sheeted ramified domain over $\C^2$.
\end{rmk}

It is unusual not to publish such an important result obtained
in the series of Papers VII--XI in 1943, which were  hand-written
but rather complete, ready for publication.
Oka probably then noticed a shadow of an unknown concept,
{\em ``Id\'eaux de domaines ind\'etermin\'es''} or {\em ``Coherence''}.
With a project in mind to settle the Levi (Hartogs' Inverse) Problem
for domains allowing singularities and ramifications,
 he would have been interested
more in inventing the new necessary notion for his project
than the publication of the important result which was enough
marvelous by itself (cf.\ \S\ref{7-11}).

As briefly mentioned at the end of \S\ref{7-11}, reading the series
of unpublished Papers VII--XI 1943 and above XII 1944,
we can see how and why Oka continued the study of the shadow of
a new notion,``{\em Coherence}'' or
 ``{\em Id\'eaux de domaines ind\'etermin\'es}''
with leaving the papers unpublished, and what he really wanted
to do; the problem of ramified Riemann domains
 left by Oka has not been settled,
although the ramification case was countered by
example (cf.\ Remark \ref{infram}).
In this sense, I
 think, the value of the series of the unpublished papers
in 1943 has not changed.
 (Here, we may recall H. Cartan's words quoted in p.\,\pageref{carw}.)

K. Oka wrote his intension 
 implicitly
 in a paragraph of Oka \cite{oka1} (\cite{oka2}) VII,
Introduction, which was written and published in an interval of
 six or eight years after Oka \cite{oka} VI 1942, and explicitly
in Oka \cite{oka} VIII, Introduction, and IX, Introduction 2 and \S23
 (cf.\ also Remark \ref{44/45}). We recall the first:
\begin{quote}\sf
Or, nous, devant le beau syst\`eme de probl\`emes \`a F. Hartogs et aux 
successeurs, voulons l\'eguer des nouveaux probl\`emes \`a ceux qui nous
suivront; or, comme le champ de fonctions analytiques de plusieurs
variables  s'\'etend heureusement aux divers branches de
math\'ematiques, nous serons permis de r\^ever divers types de nouveaux
probl\`emes y pr\'eparant.
\end{quote}
In English (from \cite{oka2} VII, translation by R. Narasimhan):
\begin{quote}\sf
Having found ourselves face to face with the beautiful problems
introduced by F. Hartogs and his successors, we should like, in turn,
to bequeath new problems to those who will follow us.
The field of analytic functions of several variables happily extends
 into divers branches of mathematics, and we might be permitted to
dream of the many types of new problems in store for us.
\end{quote}

\begin{rmk}\rm The above paragraph was deleted in the published Oka
 \cite{oka} VII without notification to K. Oka in the editorial process.
K. Oka was very unsatisfied with this change of the original text,
so that he wrote  \cite{okap} (cf.\ \cite{nog16} Chap.\ 9 `On Coherence').
\end{rmk}



The series of published papers Oka \cite{oka}, I---IX will be
classified into two groups:
\begin{itemize}
\item[(A)] I---VI+IX,
\item[(B)] VII--VIII.
\end{itemize}
In the first group he solved the Three Big Problems of Behnke--Thullen
(\S\ref{3bp}).  It is now known that 
for the solutions of those problems (even for unramified Riemann domains)
 one needs only a
rather simple {\em Weak Coherence}  (\cite{nog18}, \cite{nog21}),
 not such general Coherence Theorems proved by Oka.

The second group (B) of VII--VIII was written beyond the
Three Big Problems and was explored to a foundational theory of
modern Mathematics, not only of complex analysis
 by H. Cartan, J.-P. Serre, H. Grauert, ....

As mentioned above, the Levi (Hartogs' Inverse) Problem
for ramified domains over $\C^n$
was countered by example due to Forn{\ae}ss \cite{forn} in 1978;
in the same year K. Oka passed away.
But it is unknown the cause of the failure or what is
the sufficient condition for the validity of the problem
in ramified case;
a certain sufficient condition was lately
obtained by \cite{nog17}.

Therefore there still remains the following interesting problem:

{\bf Oka's Problem (Dream).}
   {\em What are the sufficient and/or necessary
conditions with which a ramified pseudoconvex domain over $\C^n$
is Stein?}

\begin{rmk}\rm
The English word `Dream' is taken from the above Narasimhan's
 translation of the original `r\^ever'.
The problem of the pseudoconvexity was mentioned as
the main  problem of his research from the beginning,
Oka I (\cite{oka}, published, 1936).
After settling the problem of pseudoconvexity for unramified
domains in VII---XI (unpublished, 1943) as described above,
the problem 
with interior ramification points had been mentioned
as the principal motivation  in a number of places
such as, e.g., 
 in XIII (cited above, unpublished, around 1945, cf.\ Remark
 \ref{44/45} (ii));
 in VII (\cite{oka}, published, 1948\footnote{\, It is common to use the
 year 1948 to refer Oka VII; the actual publication year is 1950
(cf.\ \cite{nog16} Chap.\ 9 `On  Coherence').}),
 the 3rd paragraph of the Introduction;
 in VIII (\cite{oka}, published, 1951),
 Introduction\footnote{\, The first two paragraphs of the Introduction
are the rephrasing of those in XIII (cited above, unpublished)
quoted in Remark \ref{44/45} (ii).} (cf.\ \S\ref{7-11} too);
 in IX (\cite{oka}, published, 1953),
 Introduction 2 and \S23
 (as already mentioned in \S\ref{7-11} above);
 and later in his lecture (\S5) at Yukawa Institute for
 Theoretical Physics, Kyoto University in March 1964 (\cite{okap}).
In particular, it is noted that the above two (published) papers VII and VIII
 are devoted to establish his {\em Three Coherence Theorems} ($\O_{\C^n}$, 
ideal sheaves of analytic subsets\footnote{\, H. Cartan gave another
 proof
for the second (cf.\ \cite{nog16} Chap.\ 9 `On  Coherence').},
 and the normalizations of
the structure sheaves of complex spaces)
;  cf.\ \cite{nog16} Chap.\ 9 `On  Coherence' and \cite{nog19}
  for more details.
\end{rmk}

Throughout his Mathematical research life,
 K. Oka seemed to intend to prove Hartogs's Inverse
 Problem (Levi's Problem) unconditionally for ramified domains over $\C^n$.
H. Grauert mentioned at his talk at
 Conference OKA 100, Kyoto/Nara 2001\footnote{\, Complex Analysis in
Several Variables---Memorial Conference of Kiyoshi Oka's Centennial
 Birthday, Kyoto/Nara 2001, Ed.\ K. Miyajima et al.,
 Advanced Studies in Pure Math.\ Vol.\ 42, Math.\ Soc.\ Jpn., Tokyo
2004. H. Grauert gave a talk titled `A simple way to perform the
 Levi-Oka-Theory', but unfortunately, his manuscript is not included
in the proceedings.}
that the Levi problem is still open for domains ramified over $\C^n$.


{\em Acknowledgment.} The author is sincerely grateful to
Mr. H. Oka for the kind agreement of the English translation
of the unpublished paper XI of \cite{okap} as the copyright holder,
and to ``Oka Kiyoshi Collection, Library of Nara Women's University''
for the resources.

%
%


\newpage
\thispagestyle{empty}
\hfill

\vspace{50mm}
\centerline{\Large\bf Part II}

\vspace{10mm}

\begin{center}
The English translation  of K. Oka's unpublished
Paper XI 1943\\
translated  by J. Noguchi.
\end{center}

\newpage
\thispagestyle{empty}
\hfill
\vspace{-10pt}
\begin{center}
\includegraphics{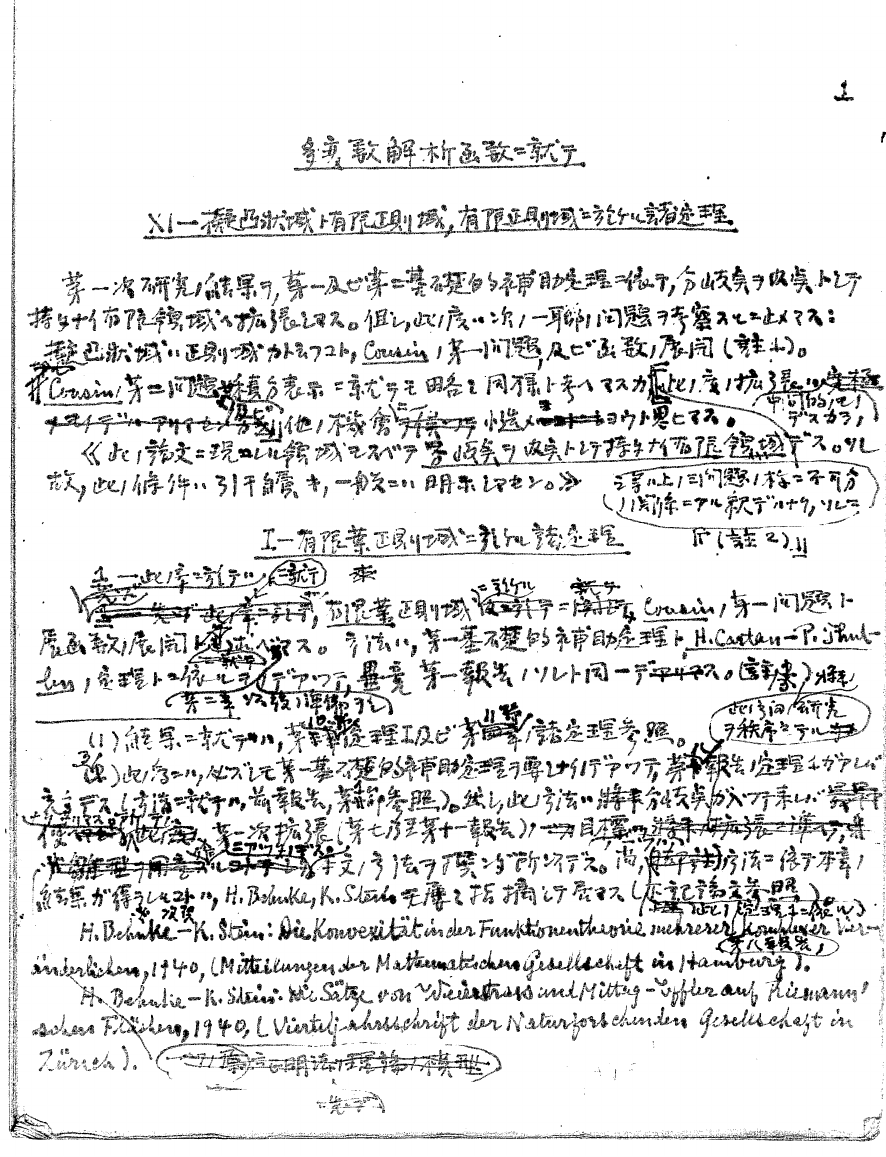}
\end{center}
\begin{center}{\small
The first page of the draft-manuscript XI, Ref.~No.~177. \\
By the courtesy of OKA Kiyoshi Collection.
Nara Women's University Library,\\
Copyright (c)\,1999; All rights reserved.
}
\end{center}

\setcounter{footnote}{0}
\newpage
\thispagestyle{empty}
\begin{center}
{\textbf{\LARGE On Analytic Functions of Several Variables}}\\
\medskip
{\Large \rm{\bf{XI}} \bf{\Large --- 
Pseudoconvex Domains and 
 Finite Domains of Holomorphy,
}}\\
{\textbf{\Large Some Theorems on Finite Domains of
Holomorphy\footnote{\,This is an English translation of the original 
Japanese text
in OKA Kiyoshi Collection, Nara Women's University Library,
Unpublished manuscript,
 http://www.lib.nara-wu.ac.jp/oka/fram/mi.html. The handwritten
 original text is found in the same Collection,
 http://www.lib.nara-wu.ac.jp/oka/moku/html/174/001.html.\\
\copyright~by courtesy of Mr.\ Hiroya Oka and OKA Kiyoshi Collection
at Nara Women's University Academic Information Center.
}
}}\\
\end{center}

\centerline{\Large Kiyoshi Oka}

\bigskip
We extend the results of the first research project
to unramified finite domains\footnote{\,Added in translation:
 Oka used the term ``finite domain'' in the sense
that it spreads over $\C^n$, not over a space with
infinity such as $\pnc$ or $(\ponec)^n$.
}
by making use of
 the First and the Second Fundamental Lemmata.
Here we restrict ourselves to deal with the following problems:
the Problem of pseudoconvex domains being domains of holomorphy,
{\em Cousin} I Problem, and Expansions of functions.\footnote{\,Cf.\
 Theorem I in \S10 and Theorems in \S11 for the results.
}

 As for {\em Cousin} II Problem and the integral
representation, we think that they will be similarly 
dealt with.\footnote{\,Since these are not in an
 inseparable relation as in the above three
theorems, and the present extension is at an intermediate stage,
we will confirm them in the next occasion.
}

In the present paper, ``{\em domains}'' are assumed to be finite and
to carry no ramification point in its interior: This assumption will be kept
all through the paper, and will not be mentioned henceforth in general.

\bigskip
\begin{center}
{\bf {I} -- Theorems in Finitely Sheeted Domains of Holomorphy}
\end{center}

{\bf \S1.} 
The present chapter describes {\em Cousin} I Problem and Expansions
of functions on finitely sheeted domains of holomorphy
 for the  preparation of what will follow in  Chapter II and henceforth.
The methods are due to the First Fundamental Lemma and
the {\em H. Cartan--P. Thullen} Theorem, and so they are 
essentially the same as those in  Report I.\footnote{\,For this
 aim the First Fundamental Lemma is not necessarily needed,
and Theorem 1 in  Report VIII suffices
(as for the methods, see \S1 of the previous Report).
This method, however, will not be effective if once a ramification point
is allowed.
Here it is noticed that one of the purposes of this first extension
(from Reports VII--XI) is to organize the
studies of this direction in future.
Because of this reason we here choose the method of the present paper.
And, it is was often mentioned also by H. Behnke and K. Stein
 that the results of the present chapter
can be obtained by the method of Theorem 1 of Report VIII
 (cf.\ the papers below).

H. Behnke--K. Stein: Approximation analytischer Funktionen in vorgegebenen 
Bereichen des Raumes von $n$ komplexen Ver\"anderlichen, 1939 (Nachrichten 
von der Gesellschaft der Wissenchaften zu G\"ottingen). 

H. Behnke--K. Stein: Die Konvexit\"at in der Funktionentheorie mehrerer 
komplexer Ver\"anderlichen, 1940 (Mitteilungen der Mathematischen 
Gesellschaft in Hamburg).  

H. Behnke--K. Stein: Die S\"atze von Weierstrass und Mittag-Leffler auf 
Riemannschen Fl\"achen, 1940 (Vierteljahrsschaft der Naturforschenden 
Gesellschaft in Z\"urich).}

We first modify (the fundamental) Lemma I to a form suitable
for our purpose. We recall it (Report VIII): 

\medskip
{\bf Lemma I.} {\em  Let $(X)$ be a 
univalent  cylinder domain in  $(x)$-space, and let $\Sig$ be
an analytic subset of $(X)$.
Let $V$ be a univalent open subset of $(X)$ with $V \supset \Sig$.
Assume that there are holomorphic functions
 $f_1 (x), f_2 (x), \ldots, f_p (x)$ in $V$ with
\[
 \Sig=\{f_1= \cdots = f_p=0\}.
\]
Let $(X^0) \Subset (X)$ be a relatively compact cylinder subdomain
and set $\Sig_0=\Sig \cap (X^0)$.

Then, for a bounded holomorphic function $\varphi (x)$ in $V$
such that $ | \varphi (x) | < M$ in $V$,
 there is a holomorphic function
$\Phi (x)$ in $(X^0)$ such that at every point of  
$\Sig_0$
$$
\Phi (x) \equiv \varphi (x) \quad (\mathrm{mod}.\ f_1, f_2, \ldots, f_p)
$$
and
$$
|\Phi (x) | < K M
$$
 on $(X^0)$. Here $K$ is a positive constant independent
from  $\varphi (x)$. }

\medskip
Let $R$ be a domain in the space of 
$n$ complex variables 
$x_1, x_2, \ldots, x_n$ (without ramification point in the interior,
 and finite) or a countable union of mutually disjoint such domains.
We consider an analytic polyhedron (a point set) $\Delta$ in $R$
satisfying the {\em following three conditions}:
\begin{enumerate}
\item[$1^\circ$]
$\Delta \Subset R$. (Therefore, $\Delta$ is contained
in a finite union of connected components of $R$,
bounded and finitely sheeted.)
\item[$2^\circ$]
$\Delta$ is defined as follows:
$$
P \!\in\! R, \ x_i \!\in\! X_i, \ f_j (P) \!\in\! Y_j \ \  
(i= 1, 2, \ldots, n; j= 1, 2, \ldots, \nu), \leqno{\quad (\Delta)}
$$
where $(x)$ is the coordinate system of the point $P$,
$X_i$ and $Y_j$ are {\em univalent} domains of
(finite) planes, and $f_j (P)$ are holomorphic functions in $R$
(in the sense of {\em one-valued} analytic functions in every
connected component of $R$; same in what follows).
\item[$3^\circ$]
The vectors $[x_1, x_2, \ldots, x_n, f_1 (P), f_2(P), \ldots, f_\nu
 (P)]$
have {\em distinct values} for {\em distinct points} of $\Delta$.
\end{enumerate}

\medskip
We introduce new variables, $y_1, y_2, \ldots, y_\nu$ 
and consider  $(x, y)$-space.
We then consider a cylinder domain,
$(X,  Y)$ with $x_i \!\in\! X_i, \, y_j \!\in\! Y_j$ 
($i= 1, 2, \ldots, n$; $j = 1, 2, \ldots, \nu$)
 together with an analytic subset
$$
y_j= f_j (P), \ P \!\in\! \Delta \qquad (j= 1, 2, \ldots, \nu).
\leqno{\quad (\Sig)}
$$
We map a point $P$ of $\Delta$ with coordinate $(x)$ 
to a point $M$ of $\Sig$ with coordinate $[x, f (P)]$.
By Condition $3^\circ$ distinct two points 
 $P_1, P_2$ of $\Delta$ are mapped always to 
distinct two points $M_1, M_2$ of $\Sig$,
and hence the map is injective.
All points of $\Sig$ is contained in 
$(X,  Y)$ and its boundary points are all
lying on the boundary of $(X, Y)$.
 (If $f_j (P) \ (j= 1, 2, \ldots, \nu)$ are simply assumed
to be holomorphic functions in $\Delta$, 
then the first half holds, but not the second half.)
Let $X^{0}_{i}, Y^{0}_{j} \ (i= 1, 2, \ldots, n; j= 1, 2, \ldots, \nu)$
be domains of complex plane such that
$X^{0}_{i} \Subset X_i, \, Y^{0}_{j} \Subset Y_j$,
and let  $\Delta_0$ denote the corresponding part of $\Delta$. 
Then, $\Delta_0 \Subset \Delta$. 
Now, let $\Delta_0 \Subset \Delta$ be an arbitrary subset.
If $P_1, P_2$ both belong to 
$\Delta_0$ and have the same coordinate,
then the distance between  $M_1, M_2$ 
carries a lower bound away from $0$. 

Let $\varphi (P)$ be an arbitrary holomorphic function in 
$\Delta$. With a point $P$ of $\Delta$ 
mapped to a point $M$ of $\Sig$, 
we consider a function $\varphi (M)$ on $\Sig$  by setting
$$
\varphi (M) = \varphi (P).
$$
As seen above, we may think a holomorphic function in $(x,y)$
defined in a univalent open set containing
$\Sig$, which agrees with 
$\varphi (M)$ on $\Sig$, 
and locally independent from $(y)$.
Therefore, Lemma I is modified to the following form:
 
\medskip
{\bf Lemma $\mathbf{I}'$.} {\em 
 Let the notation be as above.
Let $(X^{0}, Y^{0})$  be a cylinder domain 
such that $(X^{0}, Y^{0}) \Subset (X, Y)$.
Then, for a given bounded holomorphic function
$\varphi (P)$ on $\Delta$,
we may find  a holomorphic function $\Phi (x, y)$ in
$(X^{0},Y^{0})$  so that 
if $| \varphi (P)| < N$ in $\Delta$, 
$| \Phi (x, y)| < K N$  in $(X^{0}, Y^{0})$,
and  $\Phi (x, f (P))= \varphi (P)$ for all
$[x, f(P)] \in (X^{0},Y^{0}) \cap \Sig$ with coordinate $(x)$ of $P$.
Here, $K$ is a positive constant independent from 
$\varphi (P)$. }
\vspace{1mm}

\medskip
We have the following relation  between the analytic polyhedron $\Delta$ above
and a finitely sheeted domain which is convex
with respect to a family of holomorphic functions%
\footnote{\,Cf.\ the previous Report,
 \S1 for the definition of the convexity.
}: 

\medskip
{\bf Lemma 1.} {\em  
Let $\mathfrak{D}$ be a domain of holomorphy in  $(x)$-space, 
and let $\mathfrak{D}_0$ be a finitely sheeted open subset of
$\mathfrak{D}$, 
which is holomorphically convex with respect to the set
of all holomorphic functions in $\mathfrak{D}$.
For any subset $E \Subset \mathfrak{D}_0$, there exist an analytic polyhedron
 $\Delta$ and  an open subset $R$ of $\mathfrak{D}_0$ such that $E \Subset \Delta$ and
$R$ satisfies the above three Conditions,
where $f_j (P) \ (j= 1, 2, \ldots, \nu)$ 
may be taken as holomorphic functions in $\mathfrak{D}$,
 $X_i \ (i= 1, 2, \ldots, n)$  taken as disks 
$| x_i| < r$, and
$Y_j$ taken as unit disks $| y_j| < 1$.%
\footnote{\,Cf.\ the first two of the three papers of H. Behnke--K. Stein 
cited above. }
}

\medskip
{\it Proof.}~ Let $F$ be an arbitrary subset of $\mathfrak{D}_0$
 which is bounded with respect to 
$\mathfrak{D}_0$.
Since $\mathfrak{D}_0$ is finitely sheeted, 
it is immediate that
$$
F \Subset \mathfrak{D}_0 .
$$
Conversely, if $F \Subset \mathfrak{D}_0$, then 
$F$ is bounded with respect to 
$\mathfrak{D}_0$ (even if $\mathfrak{D}_0$ is not finitely sheeted).
 
Therefore, these two notions agree with each other.

We denote by 
 $(\mathfrak{F})$ the family of all holomorphic functions
 in $\mathfrak{D}$. 
Then, $\mathfrak{D}_0$ is convex with respect to $(\mathfrak{F})$, and 
$E \Subset \mathfrak{D}_0$.
As seen as above, we may take an open set $\mathfrak{D}^{'}_{0}$
with $E \subset \mathfrak{D}^{'}_{0} \Subset \mathfrak{D}_0$,
so that for every point $P_0$ of $\mathfrak{D}_0$, not belonging to
$\mathfrak{D}^{'}_{0}$, 
there is at least one function 
$\varphi (P)$ of  $(\mathfrak{F})$ satisfying 
$$
| \varphi (P_0)| > \max | \varphi (E)|.
$$
(Here, the right-hand side stands for the supremum of
$| \varphi (P)|$ on $E$.)

Let $\rho$ denote the minimum distance of $\mathfrak{D}^{'}_{0}$
with respect to $\mathfrak{D}_0$, 
and let $r$ be a positive constant such that 
any point $P (x)$ of $E$ satisfies 
$| x_i| < r \ (i=1, 2, \ldots, n)$.
We consider those points of $\mathfrak{D}_0$ such that the distance to
the boundary of $\mathfrak{D}_0$ is  $\frac{1}{2} \rho$,
and denote by $\Gamma$ the part of them over the closed polydisk $|x_i|
\leq 2r$.
As seen above, $\Gamma$ is a closed set.
It is clear that for an arbitrary point $M$ of $\Gamma$,
there are a small polydisk $(\gamma)$ with center $M$ contained in
$\mathfrak{D}$, and a function $f(P)$ of $(\mathfrak{F})$ satisfying
$$
\max |f [(\gamma)]| > 1, \quad \max | f (E)| < 1.
$$
Therefore by the Borel--Lebesgue Lemma,  $\Gamma$ is covered by finitely many
such $(\gamma)$. 
Let $f_1 (P)$, $f_2 (P)$, $\ldots$, $f_\lambda (P)$ be those
 functions associated  with them. 
Set $R = \mathfrak{D}^{(\frac {\rho}{2})}_{0}$ 
 (the set of points of $\mathfrak{D}_0$ whose distance to the boundary
of $\mathfrak{D}_0$ is greater than $\frac {\rho}{2}$.
We consider the following analytic polyhedron $\Delta$: 
$$
P \! \in \! R, \ | x_i| \!<\! r, \ | f_j (P)| \!<\! 1 
\ \ (i\!=\! 1, 2, \ldots, n; j\!=\! 1, 2, \ldots, \lambda).
 \leqno{\quad (\Delta)}
$$
Clearly, $E \subset \Delta$ and $\Delta \Subset R$. 
(The condition of Lemma requires $E \Subset \Delta$, but this is the same.)

We check Condition $3^\circ$.
Since $\mathfrak{D}$ is a domain of holomorphy, 
there is a holomorphic function whose domain of existence is
$\mathfrak{D}$. 
Let  $F (P)$ be such  one. 
Then, by the definition of domain of holomorphy%
\footnote{\,Cf.\ Behnke--Thullen's Monograph, p.\,16. 
},
for 
{mutually overlapped (the coordinates are the same)
 two points $P_1$ and $P_2$}\footnote{\,(Note by the translator.)
This means that $P_1$ and $P_2$ are distinct and
 their projections to $\C^n$ (the base points) are the same.}
 of  $\mathfrak{D}$,
the elements\footnote{\,(Note by the translator.)
That is, function elements or germs
in the present terms.}  of $F (P)$ at $P_1$ and $P_2$ 
are necessarily different.
Therefore, there exists a partial derivative
of $F(P)$ with respect to 
$x_i \ (i= 1, 2, \ldots, n)$ which takes distinct values
at $P_1$ and $P_2$, and 
the partial derivative is necessarily a holomorphic function
in  $\mathfrak{D}$.
Let $\bar{\Delta}$ denote the union of $\Delta$ and its boundary.
Since $\Delta \Subset \mathfrak{D}_0$, 
$\bar{\Delta}$ is a closed set.
Hence by the Borel--Lebesgue Lemma, there are finitely many holomorphic
functions consisting of $F (P)$ and its partial derivatives,
$$
\varphi_1 (P), \ \varphi_2 (P), \ldots, \ \varphi_\mu (P)
$$
such that the vector-valued function
$[\varphi_1 (P), \varphi_2 (P), \ldots, \varphi_\mu (P)]$
takes distinct vector-values at any two distinct points of $\bar{\Delta}$.
These functions are bounded in $\Delta$.   We set
$$
\max | \varphi_k (\Delta)| < N, \quad f_{\lambda + k} (P) = 
\frac {1}{N} \varphi_k (P) \qquad (k= 1, 2, \ldots, \mu).
$$
Then we see that the set of points of $\mathfrak{D}$
satisfying three conditions, 
$P \!\in\! R, \ | x_i| \!<\! r, \ | f_j (P)| \!<\! 1 
\ \ (i \!=\! 1, 2, \ldots, n;\,j \!=\! 1, 2, \ldots, \nu;\,\nu \!=\!
\lambda + \mu)$ 
agrees with $\Delta$.
The expression of $\Delta$ of this type satisfies
all Conditions  $1^\circ$, $2^\circ$ and $3^\circ$.
\hfill C.Q.F.D. 

Recall  that a domain of holomorphy carries the
following property:

\medskip
{\bf The First Theorem of  H.\:Cartan--P.\:Thullen.} 
{\em A finite domain of holomorphy is convex with respect to
the whole of functions holomorphic there.
}

\medskip
This theorem is an immediate consequence of the
{\em Fundamental Theorem of  H. Cartan--P. Thullen%
\footnote{\,Cf.\ Behnke--Thullen's Monograph, Chap.\ 6, \S1
and the following paper by
H. Cartan--P. Thullen : Regularit\"ats--und Konvergenzbereiche, 1932 
(Math. Annalen). } on the simultaneous
analytic continuation.}%
\footnote{\,In this way we use the Fundamental Theorem of Cartan--Thullen.
However, this theorem no longer holds if ramification points
or points of infinity are allowed to come in.
Therefore there remains a problem how to deal with these difficulties in future,
but in the present paper this theorem is not necessarily
needed in fact; cf.\ the footnote of Theorem I.
Although there do not arise no other problems of this kind,
the author  thinks that the one mentioned above is the most noticeable.
}

\bigskip
{\bf \S2.} 
We study the expansions of functions.\footnote{\,Cf.\ Report I, \S4.
}

We consider $\Delta$ in Lemma 1:
 Here we also assume that $\Delta$ satisfies the conditions
added at the end of the lemma. 
Then, $\Delta$ is of the form: 
$$
P \!\in\! R, \ | x_i| \!<\! r, \ | f_j (P)| \!<\! 1 \quad 
(i\!=\! 1, 2, \ldots, n; j\!=\! 1, 2, \ldots, \nu)  \leqno{\quad (\Delta)}
$$
We introduce complex variables 
$y_1, y_2, \ldots, y_\nu$ and in 
$(x, y)$-space we consider a polydisk 
$$
| x_i| < r, \quad | y_j| < 1 \qquad
(i= 1, 2, \ldots, n; j= 1, 2, \ldots, \nu)  \leqno{\quad (C)}
$$
and an analytic subset defined by 
$$
y_j= f_j (P), \quad P \in \Delta \qquad (j= 1, 2, \ldots, \nu) .
\leqno{\quad (\Sig)}
$$
Let $r_0$ and $\rho_0 $ be positive numbers with 
$r_0 < r$ and $\rho_0 < 1$,  
and let $\Delta_0, (C_0), \Sig_0$ respectively denote
those defined as $\Delta, (C), \Sig$
with $(r, 1)$ replaced by $(r_0, \rho_0)$.

Let $\varphi (P)$ be an arbitrary holomorphic function 
in $\Delta$. 
By Lemma $\rm{I'}$ one can construct a holomorphic function
 $\Phi (x, y)$ in $(C_0)$ such that 
 $\Phi(x, f (P))=\varphi (P)$ for all $[x, P] \in \Sig_0$. 
We expand this $\Phi (x, y)$ to a Taylor series with center
at the origin of $(C_0)$.
Then the convergence is locally uniform at every
point of $(C_0)$. 
With substituting $y_j= f_j (P) \ (j= 1, 2, \ldots, \nu)$
in that expansion, 
we obtain an expansion of $\varphi (P)$ 
in $\Delta_0$, 
whose terms are all holomorphic functions in $\mathfrak{D}$;
the convergence is locally uniform at every point of $\Delta_0$.

Since $\mathfrak{D}_0$ is the limit of
the monotone increasing sequence of 
subsets of $\mathfrak{D}_0$ satisfying the same property
as $\Delta$, 
we have the following theorem:

\medskip
{\bf Theorem 1.}~ {\em
Let $\mathfrak{D}$ be a domain of holomorphy in $(x)$-space, 
and let $\mathfrak{D}_0$ be an open subset of $\mathfrak{D}$
which is finitely sheeted and convex with respect to the whole family
$(\mathfrak{F})$ of holomorphic functions in $\mathfrak{D}$.
Then, every holomorphic function in $\mathfrak{D}_0$
is expanded to a series of functions of $(\mathfrak{F})$,
which converges locally uniformly at every point of $\mathfrak{D}_0$.
} 

\bigskip
{\bf \S3.} We next discuss {\em Cousin} I Problem.\footnote{\,Cf.\
 Report I, \S3, and the proof of Theorem I in \S5 in  Report I.}
We begin with a lemma.

\medskip
{\bf Lemma 2.} {\em 
Let $\Delta$ be as in Lemma $I'$, let $L$ be
a real hyperplane passing through a base point of $\Delta$,
 and let $S$ denote the part of $\Delta$ over $L$.  
Let $\Delta_0 \Subset \Delta$ be an open subset
and let 
$\Delta^{'}_{0}$ be the part of $\Delta_0$ in one side of $L$
and let $\Delta^{''}_{0}$ be the one in another side.
Then, for a given function  $\varphi (P)$ holomorphic 
in a neighborhood of $S$ in $R$,
one can find a holomorphic function  $\varphi_1 (P)$
(resp.\ $\varphi_2 (P)$) in $\Delta^{'}_{0}$ 
(resp.\ $\Delta^{''}_{0}$)
such that both are also holomorphic at every point
of $S$ in $\Delta_0$, and there satisfy identically
$$
\varphi_1 (P) - \varphi_2 (P) = \varphi (P).
$$}

{\em Proof.} ~ We write $x_1 = \xi + i \,\eta$ with real and imaginary parts
($i$ for the imaginary unit) 
and may assume that $L$ is defined by 
$$
\xi = 0. \leqno{\quad (L)}
$$
For $L$ is reduced to the above form by a linear transform of
$(x)$. Recall $\Delta$ to be of the following form: 
$$
P \!\in\! R, \ x_j \!\in\! X_j, \ f_k (P) \!\in\! Y_k \ \  
(j\!=\! 1, 2, \ldots, n; k\!=\! 1, 2, \ldots, \nu).
 \leqno{\quad (\Delta)}
$$
Associated with this we consider the cylinder domain $(X,Y)$ in $(x,y)$-space
as done repeatedly in above, and the  analytic subset
$\Sig$. 
Let $X^{0}_{j}, X^{1}_{j}, Y^{0}_{k}, Y^{1}_{k}$ 
be domains in the plane such that
$$
X^{0}_{j} \!\Subset\! X^{1}_{j} \!\Subset\! X_j,
 \ Y^{0}_{k} \!\Subset\! Y^{1}_{k} \!\Subset\! Y_k 
\ \  (j\!=\! 1, 2, \ldots, n; k\!=\! 1, 2, \ldots, \nu).
$$
Let $\Delta_0$ be the part of $\Delta$, where 
$(X, Y)$ is replaced by $(X^0, Y^0)$. 	      
Then, one may assume $\Delta_0$ in the lemma to be of this form.

Let $A$ be an open subset of $X_1$ in $x_1$-plane
which contains the part of the line  $\xi = 0$ in $X_1$.  
Here we take $A$ sufficiently close to this line so that
$\varphi (P)$ is holomorphic in the part of $\Delta$
over $x_1 \in A$.
Let $A_1 \Subset A$ be an open subset
which is in the same relation with respect to $X^{1}_{1}$
as $A$ to $X_1$.

By Lemma 
$\rm{I'}$ there is a holomorphic function $\Phi (x, y)$
in the cylinder domain 
with $x_1 \in A_1$ and $(x, y) \in (X^1, Y^1)$, which takes
the value  $\varphi (P)$ at every point 
$[x, f (P)]$ of $\Sig$ in this cylinder domain.
Taking a line segment or a finite union of them (closed set)
$l$ in the imaginary axis of $x_1$-plane, contained in $A_1$ 
and containing the part of the imaginary axis
inside  $X^{0}_{1}$, 
we consider {\em Cousin's integral}
$$
\Psi (x, y) = \frac {1}{2 \pi i} \int_{l} 
\frac {\Phi (t, x_2, \ldots, x_n, y)}{t - x_1}\,dt.
$$
Here the left part 
$(\xi < 0)$ of $L$ in $\Delta_0$ is 
denoted by $\Delta^{'}_{0}$,  the right part
by $\Delta^{''}_{0}$, and the orientation of the integration is 
 the positive direction of the imaginary axis.
Let $(C')$ be the part $\xi < 0$ of  $(X^0, Y^0)$,
and let $(C'')$ be that of $\xi > 0$.
Then, $\Psi (x, y)$ is holomorphic in $(C')$ 
and in $(C'')$. We distinguish  $\Psi$ as $\Psi_1$  in  $(C')$ and
that as  $\Psi_2$ in $(C'')$.
Then both of $\Psi_1$ and $\Psi_2$ are holomorphic also at
every point of $\xi = 0$ inside $(X^0, Y^0)$,
and satisfy the following relation:
$$
\Psi_1 (x, y) - \Psi_2 (x, y) = \Phi (x, y).
$$
Therefore, we obtain the required functions
$$
\varphi_1 (P) =\Psi_1 [x, f (P)], \quad \varphi_2 (P) =\Psi_2 [x, f (P)],
$$
where $(x)$ is the coordinate of a point $P$ of $R$.
\hfill C.Q.F.D.

\medskip
Let $\mathfrak{D}$ be a domain in $(x)$-space. 
Assume that for each point $P$ of $\mathfrak{D}$
there are a polydisk $(\gamma)$ with center at $P$ in
$\mathfrak{D}$  and 
a meromorphic function $g (P)$ in $(\gamma)$, 
and that the whole of them satisfies the following congruence condition:
For every pair $(\gamma_1), (\gamma_2)$ of such $(\gamma)$  
with the non-empty intersection $(\delta)$, 
the corresponding $g_1 (P)$ and $g_2 (P)$ are 
congruent in $(\delta)$; i.e., precisely, 
  $g_1 (P) - g_2 (P)$ is holomorphic in $(\delta)$.
In this way, the poles were defined in $\mathfrak{D}$.
Then, it is  the {\em Cousin I Problem} to construct a
meromorphic function $G (P)$ in $\mathfrak{D}$
with the given poles; in other words, it is congruent to
$g (P)$ in every $(\gamma)$.

Let $\mathfrak{D}$ be a finitely sheeted domain of holomorphy.
By the First Theorem of {\em Cartan--Thullen},
$\mathfrak{D}$ is convex with respect to
the family $(\mathfrak{F})$ of all holomorphic functions
in $\mathfrak{D}$. 
Therefore we may take $\mathfrak{D}$ = $\mathfrak{D}_0$ in Lemma 1,
and hence there is a $\Delta$ in $\mathfrak{D}$ stated in the lemma.
Here, however it is convenient to take a
 {\em closed analytic polyhedron $\Delta$}
with  closed bounded domains  $X_i$ and
 $Y_j \ (i= 1, 2, \ldots, n; j= 1, 2, \ldots, \nu)$.
(Naturally, $f_j (P)$ are chosen from $(\mathfrak{F})$.) 
Thus,  $\mathfrak{D}$ is a limit of a sequence of
 closed analytic polyhedra,
$$
\Delta_1, \Delta_2, \ldots, \Delta_p, \ldots,
$$
where $\Delta_p$ are such ones as  $\Delta$ above,
and  $\Delta_p \Subset E$ with the set $E$ of
all interior points of $\Delta_{p + 1}$.

Now, we take a $\Delta_p$ and divide it into $(A)$
as stated in \S3 of the previous Report\footnote{\,(Note by the translator.)
 The ``previous Report'' is ``Report X''; there in \S3, small
closed cubes are defined so that their sides are parallel to real and imaginary
axes of the complex coordinates of the base space $\C^n$.}:
Here, we choose $2n$-dimensional closed cubes
for $(A)$ and its base domain  $(\alpha)$. 
We also allow some of $(A)$ to be of incomplete form,
and take $(A)$ sufficiently small so that $(A) \Subset (\gamma)$
for every  $(A)$ with one of $(\gamma)$ above. 
Choosing arbitrarily such  $(\gamma)$, 
we associate $g (P)$ with $(\gamma)$, and then $g(P)$ with this $(A)$. 

Let $(A)_1, (A)_2$ be a pair of $(A)$
adjoining by a face (a $(2n-1)$-dimensional closed cube). 
The meromorphic functions $g_1 (P)$ and $g_2 (P)$ associated with
them are congruent 
 in a neighborhood of the common face 
(a neighborhood in $\mathfrak{D}$, same in below).
It follows from Lemma 2 that there is a meromorphic function
with the given poles in a neighborhood of the union
$(A)_1 \cup (A)_2$.
It is the same for a union of $(A)$ such as, e.g.,
$$
\left(\alpha^{(1)}_{j, q}, \alpha^{(2)}, \ldots, \alpha^{(n)}
\right),
$$
where $\alpha$ are closed squares, $q$ and 
$\alpha^{(2)}, \ldots, \alpha^{(n)}$ are given ones, 
and $j$ is arbitrary. Here $(A)$ may be disconnected.
Repeating this procedure, we obtain a meromorphic function $G (P)$
in a neighborhood of $\Delta_p$ with the given poles.

Thus, we have
$$
G_1 (P), G_2 (P), \ldots, G_p (P), \ldots .
$$ 
We consider
$$
H (P) = G_{p + 1} (P) - G_p (P).
$$
Then, $H (P)$ is a holomorphic function in a neighborhood of
$\Delta_p$. Hence by Theorem 1, 
$H_p$ is expanded to a series of functions of $(\mathfrak{F})$
which converges uniformly in a neighborhood of $\Delta_p$.
By this we immediately see the existence of a meromorphic function
$G (P)$ in $\mathfrak{D}$ with the given poles.
(The method of the proof is exactly the same as in the case of
 univalent cylinder domains.)
Thus we obtain the following theorem.

\medskip
{\bf Theorem 2.} {\em In a finitely sheeted domain of holomorphy,
the Cousin I Problem is always solvable.}

\begin{center}
\bigskip
{\bf II --- The Main Problem}
\end{center}

{\bf \S4.} In this chapter 
 we solve the main part of the problem
 abstracted from  the series of those discussed at the beginning by virtue of
 the First Fundamental Lemma.\footnote{\,Except for
 the use of this lemma, the content is essentially the
same as in Report VI, Chap.~1.
}

We begin with explaining the problem.
Let $\mathfrak{D}$ be a {\em bounded finitely sheeted domain} 
in $(x)$-space. 
We consider a real hyperplane 
 with non-empty intersection
with the base domain of $\mathfrak{D}$.
We write $x_1$ 
as
$$
x_1= \xi + i \,\eta .
$$
For the sake of simplicity we assume that this hyperplane
is given by $\xi = 0$.
Let $a_1, a_2$ be real numbers such that
$$
a_2 < 0 < a_1,
$$
and the hyperplanes $\xi = a_1, \ \xi = a_2$ have
both non-empty intersections with the base domain of $\mathfrak{D}$. 
Let $\mathfrak{D_1}$ (resp.\ $\mathfrak{D_2}$) denote
the part of $\xi < a_1$ (resp.\ $\xi > a_2$) in $\mathfrak{D}$,
and let  $\mathfrak{D}_3$ be the part of $a_2 < \xi < a_1$
in  $\mathfrak{D}$.
We assume that {\em every connected component of $\mathfrak{D}_1$
and $\mathfrak{D}_2$ is a domain of holomorphy.}
Then, necessarily so is every component of $\mathfrak{D}_3$.

Let $f_j (P) \ (j= 1, 2, \ldots, \nu)$ be {\em holomorphic functions
in $\mathfrak{D}_3$.} 
We consider a subset $E$ of $\mathfrak{D}$
such that $E \supset \mathfrak{D}\setminus \mathfrak{D}_3$
and the following holds: A point $P$ of $\mathfrak{D}_3$ belongs to $E$
if and only if
$$
| f_j (P)| < 1 \qquad (j= 1, 2, \ldots, \nu).
$$
We assume that $E$ has  connected components which
extend over the part $\xi < a_2$ and over $\xi > a_1$. 
Let $\Delta$ be such one of them. 

{\em We assume the following three conditions for this $\Delta$:

\begin{enumerate}
\item[$1^\circ$]
Let $\delta_1$ be a real number such that 
$0 < \delta_1 < \min\{a_1,  \,-a_2\}$.
Let $A$ denote the set of point $P (x)$ of $\Delta$ with
$ | \xi| < \delta_1$. Then,
$$
A \Subset \mathfrak{D}. 
$$
\item[$2^\circ$]
Let $\delta_2$ be a positive number and let 
$\varepsilon_0$ be a positive number less than $1$.
For every $p$ of 
$1, 2, \ldots, \nu$, any point $P$ of $\mathfrak{D}_3$ satisfying
$$
| f_p (P)| \geq  1 - \varepsilon_0
$$
does not lie over
\begin{center}
$| \xi - a_1 | < \delta_2$ \quad or \quad $| \xi - a_2 | < \delta_2$.
\end{center}
\item[$3^\circ$]
The vector-values
$$
[f_1 (P), f_2 (P), \ldots, f_\nu (P)]
$$
are never identical for
{mutually overlapped two points} of $A$.
\end{enumerate}
}

By the second Condition, $\Delta$ is a domain. 
Let $\rho_0$ be a real number such that 
$1 - \varepsilon_0 < \rho_0 < 1$,
and consider a subset
 $\Delta_0$ of $\Delta$ such that  
$\Delta_0 \supset \Delta \setminus \mathfrak{D}_3$ and 
for a point of $\mathfrak{D}_3 \cap \Delta$
it belongs to $\Delta_0$ if and only if
$$
| f_j (P)| < \rho_0 \qquad ( j= 1, 2, \ldots, \nu).
$$
By Condition $2^\circ$, $\Delta_0$ is an open set. 
Denote by $\Delta^{'}_{0}$ (resp.\ $\Delta^{''}_{0}$)
the part of $\xi < 0$ (resp.\ $\xi > 0$) 
in $\Delta_0$.

The theme of the present chapter is the following problem.

\medskip
{\em
Let the notation be as above.
Let  $\varphi (P)$ be a given holomorphic function in $A$.
Then, construct  holomorphic functions, $\varphi_1 (P)$ in $\Delta^{'}_{0}$
and $\varphi_2 (P)$ in $\Delta^{''}_{0}$, which are holomorphic
in the part of $\Delta_0$ over $\xi = 0$, and identically satisfy 
$$
\varphi_1 (P) - \varphi_2 (P) = \varphi (P).
$$
}

\bigskip
{\bf \S5.} 
By making use of the method of Lemma 2 we first solve a part
of the problem related to $\mathfrak{D}_3$.
Let $y_1, y_2, \ldots, y_\nu$ be complex variables,
and consider in  $(x, y)$-space the analytic subset
$$
y_k= f_k (P), \quad P \in \mathfrak{D}_3 \qquad (k= 1, 2, \ldots, \nu)  .
\leqno{\quad (\Sig)}
$$
Let $r$ and $r_0$ be  positive numbers with 
$r_0 \!<\! r$, and let $r_0$ be taken sufficiently large so that 
the bounded domain $\mathfrak{D}$
is contained in the polydisk of radius 
$r_0$ with center at the origin. 
Let $\rho$ be a number with 
$\rho_0 < \rho < 1$, and consider  polydisks 
$$  | x_j| < r, \,| y_k| < \rho \quad (j= 1, 2, \ldots, n; 
k= 1, 2, \ldots, \nu), \leqno{\quad (C)}
$$
 and 
$$
 | x_j| < r_0, \,| y_k| < \rho_0 \quad (j= 1, 2, \ldots, n; 
k= 1, 2, \ldots, \nu).
 \leqno{\quad (C_0)} 
$$ 
Let $\delta$ be a positive number with 
$\delta < \delta_1$, and consider a set 
$$
P \in A, \quad | \xi| < \delta, \quad | f_k (P)| < \rho 
\qquad (k= 1, 2, \ldots, \nu).  \leqno{\quad (A')}
$$

Since $\varphi (P)$ is holomorphic in 
$A$, by Lemma $\rm{I'}$ we can construct a holomorphic function 
$\Phi (x, y)$ in the intersection of $(C)$ and 
$| \xi| < \delta$ such that 
$\Phi[x, f (P)]=\varphi(P)$ 
for $[x, f(P)] \in \Sig$ with  $P \in A'$ and the coordinate $x$ of $P$.
We take a line segment $l$ (connected and closed)
in the imaginary axis of $x_1$-plane,
so that it is contained in the disk  $| x_1| < r$ and
the both ends are out of the disk $| x_1| < r_0$. 
We then consider the {\em Cousin integral}\label{coin}
$$
\Psi (x, y) = \frac {1}{2 \pi i} \int_{l} 
\frac {\Phi (t, x_2, \ldots, x_n, y_1, \ldots, y_\nu)}{t - x_1} \,dt,
$$
where the orientation is in the positive direction of the imaginary axis.

Substituting $y_k = f_k (P)$ in $\Psi (x, y)$ , we get
$$
\psi (P) = \frac {1}{2 \pi i} \int_{l} 
\frac {\Phi [t, x_2, \ldots, x_n, f_1 (P), \ldots, f_\nu (P)]}{t - x_1} \,dt.
$$
The function $\psi (P)$ represents respectively
a holomorphic function $\psi_1 (P)$ in 
$\Delta^{'}_{0} \cap \mathfrak{D}_3$ and $\psi_2 (P)$ 
in $\Delta^{''}_{0} \cap \mathfrak{D}_3$. 
These are also holomorphic at every point of $\Delta_0$ 
over $\xi = 0$, and satisfy the relation: 
$\psi_1 (P) - \psi_2 (P) = \varphi (P)$. 

We modify a little the expression of this solution.
We draw a circle $\Gamma$ of radius $\rho_0$ with center at the origin
in the complex plane.
It follows from {\em Cauchy} that 
for  $| \xi| < \delta, \, | x_j| < r$ and 
$| y_k| < \rho_0 \ (j= 1, 2, \ldots, n; k= 1, 2, \ldots, \nu)$
$$
\Phi (x, y) = \frac {1}{(2 \pi i)^\nu} \int_{\Gamma} \int_{\Gamma} 
\cdots \int_{\Gamma} 
\frac {\Phi (x_1, \ldots, x_n, u_1, \ldots, u_\nu)}{(u_1 - y_1) \cdots 
(u_\nu - y_\nu)} \,d u_1 d u_2 \cdots d u_\nu,
$$
where the integral is taken on $\Gamma$ with the positive orientation.
We write this simply as follows:
$$
\Phi (x, y) = \frac {1}{(2 \pi i)^\nu} \int_{(\Gamma)} 
\frac {\Phi (x, u)}{(u_1 - y_1) \cdots (u_\nu - y_\nu)} \,du.
$$
We substitute $y_k = f_k (P) \ (k = 1, 2, \ldots, \nu)$
in this integral expression of $\Phi (x, y)$, 
 change $x_1$ with $t$, and substitute them in the integral expression
of $\psi (P)$ above. 
Then, with $t = u_0$ we obtain 
$$
\psi (P) = \int_{(l, \Gamma)} \chi (u, P) \Phi (x', u) \,du,
 \leqno{\quad (1)}
$$
$$
\chi (u, P) = \frac {1}{(2 \pi i)^{\nu + 1} 
(u_0 - x_1) [u_1 - f_1 (P)] \cdots [u_\nu - f_\nu (P)]} .
$$
Here we simply write $\Phi (x', u)$ 
for $\Phi (u_0, x_2, \ldots, x_n, u_1, \ldots, u_\nu)$,
and use the same simplification for the integral symbol as above:
It will be clear without further explanation.
Then we can use this $(1)$ in
$\Delta_0 \cap \mathfrak{D}_3$ for the integral expression
of $\psi (P)$ above. 

\bigskip
{\bf \S6.} 
There are univalent domains of holomorphy in  $(u)$-space, which contain
  the closed cylinder set  $(l, \Gamma)$ with
 $\ u_0 \in l, \, u_k \in \Gamma \ (k = 1, 2, \ldots, \nu)$,
and are arbitrarily close to $(l, \Gamma)$.
Let $V$ be  such one of them.
We shall take $V$ sufficiently close to $(l, \Gamma)$, 
as we will explain at each step in below.

Firstly, we would like to {\em construct a meromorphic function
 $\chi_1 (u, P)$ in
 $(V, \mathfrak{D}_1) \ \ ((u) \in V, \ P (x) \in \mathfrak{D}_1)$,
with the same poles as $\chi (u, P)$ of {\rm (1)} in $(V, \mathfrak{D}_3)$
and without other poles.}

This is possible by Theorem 2, because
$(V, \mathfrak{D}_1)$ is a finitely sheeted domain of holomorphy,
and for the pole distribution
the congruent condition is satisfied with $V$ sufficiently close
to $(l, \Gamma)$ by Condition $2^\circ$ on $\Delta$.

Note that $\chi - \chi_1$ is holomorphic in $(V, \mathfrak{D}_3)$. 
By the First Theorem of 
{\em Cartan--Thullen} 
$(V, \mathfrak{D}_3)$ is convex with respect to 
the family of all holomorphic functions in $(V, \mathfrak{D}_1)$.
By Theorem 1, $\chi - \chi_1$ is hence expanded to a series 
of holomorphic functions in  $(V, \mathfrak{D}_1)$, 
convergent locally uniformly at every point of $(V, \mathfrak{D}_3)$.
Therefore,  taking  $V$ closer to $(l, \Gamma)$, 
we have the following function $F_1(u, P)$
for a positive number $\varepsilon$: 
{\em $F_1 (u, P)$ is holomorphic in 
$(V, \mathfrak{D}_1)$ and for the analytic polyhedron
$A$ given in \S4,} 
$$
| \chi - \chi_1 - F_1| < \varepsilon \quad \hbox{in } (V, A). 
$$
Put
$$
K_1 (u, P) = \chi - \chi_1 - F_1.
$$
The function $K_1 (u, P)$ is holomorphic in $(V, \mathfrak{D}_3)$, and 
$| K_1| < \varepsilon$ in $(V, A)$. 
For $\mathfrak{D}_2$, we construct $K_2 (u, P)$, similarly. 
With these preparations we change the integration (1) as follows:
$$
I_1 (P) = \int_{(l, \Gamma)} [\chi (u, P) - K_1 (u, P)]\Phi (x', u) \,du,
\leqno{\quad (2)}
$$
$$
I_2 (P) = \int_{(l, \Gamma)} [\chi (u, P) - K_2 (u, P)]\Phi (x', u) \,du.
$$

If $(u) \in (l, \Gamma)$, then 
$\chi - K_1$ is equal to $\chi_1 + F_1$, so that  
it is meromorphic in $P (x) \in \mathfrak{D}_1$,
and in particular, it is holomorphic in $\Delta^{'}_{0}$.
Therefore, {\em $I_1 (P)$ is holomorphic in 
$\Delta^{'}_{0}$}; 
similarly, {\em $I_2 (P)$ is holomorphic in $\Delta^{''}_{0}$.} 
\vspace{1mm} 

The analytic functions $I_1 (P)$ and $I_2 (P)$ are holomorphic at 
every point of $\Delta_0$ over $\xi = 0$:
For $\psi (P)$ in (1) has this property and 
the both of $K_1$ and $K_2$ are holomorphic functions. 
By the property of $\psi (P)$, 
the functions $I_1 (P)$ and $I_2 (P)$ satisfy the {\em following relation}:
$$
I_1 (P) - I_2 (P) = \varphi (P) - \int_{(l, \Gamma)} 
[K_1 (u, P) - K_2 (u, P)] \Phi (x', u) \,d u. \leqno{\quad (3)}
$$
We write
$$
K (u, P) = K_1 (u, P) - K_2 (u, P).
$$
Observing this identity again,
 we see that $\varphi (P)$ is a holomorphic function in $P \in A$, 
$K$ is a holomorphic function in $(u) \!\in\! V$ and
 $P \!\in\! \mathfrak{D}_3$, 
and $\Phi (x, y)$ is a holomorphic function in 
$(x, y) \in (C)$ with $| \xi| < \delta$.
Therefore, the right-hand side is {\em a holomorphic function
in $P(x) \in A$};
 hence, {\em it is the same for the left-hand side as above.}
Put
$$
\varphi_0 (P) = I_1 (P) - I_2 (P).
$$

Let $\varphi_0$ and $K$ be given functions, and
let  $\varphi, \Phi$ be a pair of unknown functions
satisfying the relations described next
 below\footnote{\,(Note by the translator.) Here is a point of the
arguments of the proof, but one must be careful of
the notational confusion with $\varphi$, $\varphi_0$ and $\Phi$
 discussed already.}.
We consider a {\em functional equation}
$$
\varphi (P) = \int_{(l, \Gamma)} K (u, P) \Phi (x', u) \,d u + \varphi_0
(P).
\leqno{\quad (4)}
$$
\begin{itemize}
\item[{}]
{\em Here, $\Phi (x', u)$ stands for $\Phi (u_0, x_2, \ldots, x_n, u_1, \ldots, u_\nu)$,
$\varphi_0 (P)$ is a holomorphic function in $A$, 
and $K (u, P)$ is  a holomorphic function in $(V, \mathfrak{D}_3)$.
In $(V, A)$, $| K (u, P)| < 2\varepsilon$. 
For the unknown functions $\varphi(P)$ and $\Phi (x, y)$, 
the following condition is imposed besides (4):
 $\varphi (P)$ is a holomorphic function in $P \in A$, 
$\Phi (x, y)$ is a holomorphic function in 
$(x, y) \in (C)$ with $| \xi| < \delta$, 
and for every point $[x, f (P)]$ of  $\Sig$ with $P \in A'$,
$\Phi (x, y) = \varphi (P)$.}
\end{itemize}
Since these conditions are imposed, this functional equation
is not so different from the definite integral equation.

We are going to show that this equation has necessarily a solution
for a sufficiently small  $\varepsilon$.
Before it we confirm that it suffices for our end.
Suppose that there exist  functions
 $\varphi (P)$ and $ \Phi (x, y)$ as above. 
Substitute $\Phi (x', u)$ to (2).
The function $I_1 (P)$ thus obtained 
is clearly holomorphic in $\Delta^{'}_{0}$.
Similarly, $I_2 (P)$ is holomorphic in $\Delta^{''}_{0}$. 
It is clear that these analytic functions are also holomorphic
at every point of $\Delta_0$ over $\xi= 0$.
One easily sees relation (3) among them.
(The argument above is just a repetition of a deduction once done
with clarifying the conditions.)
Thus, these $I_1(P)$ and $I_2(P)$ are the solutions of the problem
described in \S4.
As seen above, it suffices to solve equation (4); here one may take
$\varepsilon$ as small as necessary.

Now, we solve equation (4). 
Recall that the analytic polyhedron $A$ 
is of the following form:
$$
P \!\in\! \Delta, \quad | \xi| \!<\! \delta_1, \quad | f_k (P)| \!<\! 1 
\qquad (k \!=\! 1, 2, \ldots, \nu).  \leqno{\quad (A)}
$$
Moreover, the analytic polyhedron $A'$ is obtained
by replacing  $(\delta_1, 1)$ of $A$ by $(\delta, \rho)$ with
$0 \!<\! \delta \!<\! \delta_1$ and $\rho_0 \!<\! \rho \!<\! 1$.
Taking $(\delta', \rho')$ with
 $\rho \!<\! \rho' \!<\! 1$ and $\delta \!<\! \delta' \!<\! \delta_1$,
we define an analytic polyhedron $A''$, replacing 
$(\delta, \rho)$ by  this $(\delta', \rho')$ in the definition of $A'$.
We have the following relation among them:
$$
A' \Subset A'' \Subset A.
$$

The function $\varphi_0 (P)$ is holomorphic in 
$A$, and hence bounded on $A''$. 
Suppose that 
$$
| \varphi_0 (P)| < M_0 \quad \mbox{on}~ A''.
$$
We denote by $(C')$ the cylinder domain given by 
$(x, y) \!\in\! (C)$ and $| \xi| \!<\! \delta$.
By {\em Lemma $I'$} we can take a holomorphic function
$\Phi_0 (x, y)$ in $(C')$ so that 
it has  values $\varphi_0 (P)$ at 
 points $[x, f (P)]$ of $\Sig$ with $P \!\in\! A'$,
and
$$
| \Phi_0 (x, y)| < N M_0  \quad \hbox{on}~ (C'),
$$
where $N$ is a positive constant independent of $\varphi_0 (P)$ 
(also independent of 
$M_0$, and of  $\varphi_0 (P)$ being holomorphic in $A$).
 Applying the operator $K (\Phi_0)$ for $\Phi_0 (x, y)$
defined by 
$$
\varphi_1 (P) = K (\Phi_0) = \int_{(l, \Gamma)} K (u, P) \Phi_0 (x', u)
\,d u, 
$$
we construct a function $\varphi_1 (P)$. 
For $(u) \in (l, \Gamma)$, 
$K(u, P)$ is holomorphic in $P (x) \in \mathfrak{D}_3$, 
and $\Phi_0 (x', u)$ is holomorphic 
in $| x_j| < r \ (j = 2, 3, \ldots, n)$, and so is in $(C)$. 
Hence, {\em $\varphi_1 (P)$  
is a holomorphic function in $\mathfrak{D}_3$.
} 

We next estimate $\varphi_1 (P)$. 
For $(u) \!\in\! (l, \Gamma)$ and $P \!\in\! A$, 
$| K (u, P)| < 2 \varepsilon$, and 
$| \Phi_0 (x, y)| < N M_0$ in $(C')$. 
Therefore,  we have in $A$,
$$
| \varphi_1 (P)| < 2 \varepsilon N N_1 M_0, \quad N_1 = 2 r (2 \pi
\rho_0)^\nu. 
$$
Therefore in first we take $\varepsilon$ so that
$$
2 \varepsilon N N_1 = \lambda < 1.
$$

\medskip
Thus, 
{\em 
$\varphi_1 (P)$ is a bounded holomorphic function
in $A$, and necessarily so is in $A''$. 
}
As we choose a function $\Phi_0 (x, y)$ for $\varphi_0 (P)$, 
we choose a function $\Phi_1 (x, y)$ for $\varphi_1 (P)$,
and by setting  $\varphi_2 (P) = K (\Phi_1)$, we construct 
$\varphi_2 (P)$. 
Inductively, we obtain
$\varphi_p (P)$ and $\Phi_p (x, y)$ ($p=0,1,2, \ldots$).
Then we consider the following function series: 
$$
\varphi_0 (P) + \varphi_1 (P) + \cdots + \varphi_p (P) +
\cdots,\vspace{-4mm} 
\leqno{\quad (5)}
$$
$$
\Phi_0 (x, y) + \Phi_1 (x, y) + \cdots + \Phi_p (x, y) + \cdots .
\leqno{\quad (6)}
$$

It follows that $\varphi_p (P)$ is holomorphic in $\mathfrak{D}_3$, and 
$\Phi_p (x, y)$ is holomorphic in $(C')$. 
In $A$, 
$$
| \varphi_p (P)| < \lambda^p M_0 \quad (p > 0), 
$$
and in $(C')$, 
$$
| \Phi_p (x, y)| < \lambda^p N M_0.
$$
Therefore,  (5) (resp.\ (6)) converges uniformly in 
$A$ (resp.\ $(C')$).
We denote the limits by $\varphi (P)$ and $\Phi (x, y)$,
respectively. 
We see that $\varphi (P)$ (resp.\ $\Phi (x, y)$) 
is holomorphic in $A$ (resp.\ $(C')$).
Since  $\Phi_p (x, y) \ (p = 0, 1, \ldots)$ take values 
$\varphi_p (P)$ at points $[x, f (P)]$ of $\Sig$ with
$P \in A'$, $\Phi (x, y)$ there takes values $\varphi (P)$. 
Therefore, it suffices to show that
 $\varphi (P)$ and $\Phi (x', u)$ satisfy functional equation (4) 
in $P \!\in\! A$. 
Now for $P \!\in\! A$ we have 
$$
\varphi_0 \!=\! \varphi_0, \ \varphi_1 \!=\! K (\Phi_0), \ \varphi_2
\!=\! K (\Phi_1), \ldots, \varphi_{p + 1} \!=\! K (\Phi_p), \ldots, 
$$
so that
$$
\varphi = K (\Phi) + \varphi_0.
$$
{\em Thus, the problem stated at the end of \S4 is always solvable.}

\bigskip
\begin{center}
{\bf
III --- Pseudoconvex domains and domains of holomorphy, 
 theorems on domains of holomorphy
}
\end{center}

{\bf \S7.} Apart from the theme we prepare some lemmata
for a moment (\S\S7--9).

 We begin with reformulating the Second Fundamental Lemma.

\medskip
{\bf Lemma II.} 
{\em  Let $\mathfrak{D}$ be a finite unramified pseudoconvex domain
over $(x)$-space.
Then there necessarily exists a real-valued continuous function
$\varphi_0 (P)$, satisfying the following two conditions:
\begin{enumerate}
\item[$1^\circ$]
For every real number $\alpha$, 
$\mathfrak{D}_\alpha \Subset \mathfrak{D}$, where $\mathfrak{D}_\alpha$
denotes the set of all points $P \in \mathfrak{D}$ with
$\varphi_0 (P) < \alpha$.
\item[$2^\circ$] 
In a neighborhood $U$ of every point 
$P_0$ of $\mathfrak{D}$, there is a hypersurface $\Sig \subset U$,
passing through $P_0$ such that
$\varphi_0 (P) > \varphi_0 (P_0)$ for $P \in \Sig \setminus  P_0$.
\end{enumerate}
}

\medskip
{\em Proof.}  
As a consequence of the former Report we know that there is
 a pseudoconvex function in
$\mathfrak{D}$ satisfying Condition $1^\circ$ and Condition $2^\circ$
outside of an exceptional discrete subset without accumulation point
in $\mathfrak{D}$.
Let $\varphi (P)$ be a such function, and let 
$E_0$ denote the exceptional discrete subset, provided that
it exists. If there is a point of $E_0$ on
$\varphi(P) = \lambda$ for $\lambda\in \mathbf{R}$,
we then call $\lambda$ an {\em exceptional value} of  $\varphi (P)$. 
For an arbitrary real number $\alpha$, 
we denote by $\mathfrak{D}_\alpha$ the set of all points
 $P \in \mathfrak{D}$ with $\varphi (P) < \alpha$.
Since $\mathfrak{D}_\alpha \Subset \mathfrak{D}$ by  Condition $1^\circ$,
$\mathfrak{D}_\alpha$ is bounded and finitely sheeted.
This remains valid for a little bit larger $\alpha$,
and so there are only finitely many points
of $E_0$ in $\mathfrak{D}_\alpha$.
Since
 $\lim_{\alpha \rightarrow 
\infty} \mathfrak{D}_\alpha = \mathfrak{D}$,
the set of the exceptional values is countable.
Let the exceptional values be
$$
\lambda_1, \lambda_2, \ldots, \lambda_p, \ldots, \qquad \lambda_p < 
\lambda_{p + 1}.
$$

Let $\alpha_0$ be a {non-exceptional value}
and set $\mathfrak{D}_{\alpha_0} = \Delta$. 
In $\Delta$ we consider 
$$
\psi (P) = - \log d (P).
$$
Here $d (P)$ denotes the Euclidean boundary distance function
with respect to $\Delta$, 
and the logarithm symbol stands for the real branch.
Since $\Delta$ is bounded, $\psi (P)$ is a continuous function. 
For any real number $\alpha$,
we denote by $\Delta_\alpha$ the set of all points $P$ of $\Delta$
with $\psi (P) < \alpha$. 
Then, $\Delta_\alpha \Subset \Delta$. 
Thus, $\psi (P)$ satisfies Condition $1^\circ$ in $\Delta$. 
We next check Condition $2^\circ$.
Let $P_0$ 
be an arbitrary point of $\Delta$, and set $\psi (P_0) = \beta$. 
We draw a $2n$-dimensional ball  $S$ of radius $e^{- \beta}$ with center $P_0$
in $\mathfrak{D}$.
Then, $S \subset \Delta$ and there is a point $M$ on the
boundary of $S$, satisfying $\varphi (M) = \alpha_0$. 
Since $\varphi (P)$ satisfies Condition~$2^\circ$ in a neighborhood 
of $\varphi (P) = \alpha_0$, there is a complex hypersurface $\sigma$ in a 
neighborhood of $M$, passing through $M$, such that
$\varphi (P) > \alpha_0$ for $P \in \sigma \setminus \{M\}$.
By a parallel translation
$$
x'_{i} = x_i + a_i \qquad (i = 1, 2, \ldots, n),  \leqno{\quad (T)}
$$
we move $M$ to $P_0$, and $\sigma$  to $\sigma'$. 
Then, $\sigma'$ is defined in a neighborhood of $P_0$. 
Let $P'$ be a point of $\sigma'$ different to $P_0$.
Then the corresponding point $P$ of $\sigma$ 
lies in $\varphi (P) > \alpha_0$,
and the (Euclidean) distance between 
$P$ and $P'$  is $e^{- \beta}$, so that 
if $P'$ belongs to $\Delta$, $P'$ lies in the part of
  $\psi (P) > \alpha$\footnote{\,(Note by the translator.)
Here $P$ is used in a different sense from the one just before
in the same sentence, and $\alpha$ is a typo of $\beta$.
They should be read as``$\psi > \beta$''}. 
Therefore, 
{\em $\psi (P)$  is a continuous function in $\Delta$, 
satisfying Conditions $1^\circ$ and $2^\circ$.}

\medskip
We take a sequence of real numbers, 
$\mu_1, \mu_2,\ldots, \mu_p, \ldots$ such that 
$$
\mu_1 < \lambda_1, \quad \lambda_p < \mu_{p + 1} < \lambda_{p + 1}.
$$
Taking $\alpha_0$ with 
$$
\lambda_1 < \alpha_0 < \mu_2,
$$
we consider  $\psi (P)$ above. 
Choosing $\alpha_0$ sufficiently close to $\lambda_1$, 
we may take $\beta$ for this $\psi (P)$, satisfying 
$$
\mathfrak{D}_{\mu_1} \Subset \Delta_{\beta} \Subset \mathfrak{D}_{\lambda_1}.
$$
Modifying $\varphi (P)$ by making use of $\psi (P)$
thus obtained (similarly to the last part of the previous Report),
we construct $\varphi_1 (P)$:  We explain it in below.

Let $\beta_1, \beta_2$ be real numbers with the same property as  
$\beta$ above such that 
$$
\beta_1 < \beta_2 .
$$
Let $\gamma_1, \gamma_2$ be real numbers with 
$$
\lambda_1 < \gamma_1 < \gamma_2 < \alpha_0.
$$
We divide $\mathfrak{D}$ into five parts 
$\mathfrak{D}_j$ ($j = 1, 2, \ldots, 5$) defined by 
$$
\mathfrak{D}_1 = \Delta_{\beta_1}, \quad 
\Delta_{\beta_1} \cup \mathfrak{D}_2 = \Delta_{\beta_2}, \quad 
\Delta_{\beta_2} \cup \mathfrak{D}_3 = \mathfrak{D}_{\gamma_1},  
$$
$$
\mathfrak{D}_{\gamma_1} \cup \mathfrak{D}_4 = \mathfrak{D}_{\gamma_2}, \quad 
\mathfrak{D}_{\gamma_2} \cup \mathfrak{D}_5 = \mathfrak{D}.
$$
By taking a suitable $B$ and a sufficiently large
positive $A$, 
we have
$$
\Psi (P) = A [\psi (P) - B]
$$
satisfying 
\begin{align*} 
\varphi (P) &> \Psi (P) \qquad \hbox{in } \mathfrak{D}_1, \\
\varphi (P) &< \Psi (P) \qquad \hbox{in } \mathfrak{D}_3 \cup
 \mathfrak{D}_4.
\end{align*}
\noindent 
Also by taking a suitable real number $B'$ and a 
sufficiently large positive number $A'$, 
we have
$$
\Phi (P) = A' [\psi (P) - B']
$$
satisfying 
\begin{align*}
\Psi (P) &> \Phi (P) \qquad \hbox{in } \mathfrak{D}_3, \\
\Psi (P) &< \Phi (P) \qquad \hbox{in } \mathfrak{D}^{'}_{5}, \\
\varphi (P) &< \Phi (P) \qquad \hbox{in } \mathfrak{D}_5,
\end{align*}
where $\mathfrak{D}^{'}_{5}$ is the  part of 
$\mathfrak{D}_5$ (a neighborhood)  containing the point set, 
$\varphi (P) \!=\! \gamma_2$.  
We define $\varphi_1 (P)$ as follows: 
\[
\begin{array}{rclll}
\varphi_1 (P) &=& \varphi (P) & \hbox{in} & \mathfrak{D}_1, \\
\varphi_1 (P) &=& \max [\varphi (P), \Psi (P)] &
 \hbox{in} & \mathfrak{D}_2, \\
\varphi_1 (P) &=& \Psi (P) & \hbox{in} & \mathfrak{D}_3, \\
\varphi_1 (P) &=& \max [\Psi (P), \Phi (P)] & \hbox{in} & \mathfrak{D}_4,\\
\varphi_1 (P) &=& \Phi (P) & \hbox{in} & \mathfrak{D}_5.
\end{array}
\]

We examine $\varphi_1 (P)$ thus defined.
It follows that $\varphi_1 (P)$ 
is a real one-valued function in $\mathfrak{D}$,
which is clearly continuous.
Since $\psi (P)$ satisfies Condition~$2^\circ$, and 
$\varphi (P)$ satisfies Condition~$2^\circ$ outside a set 
of exceptional points without accumulation point in
$\mathfrak{D}$, 
$\varphi_1 (P)$ satisfies the same condition as 
$\varphi (P)$. 
We check up the exceptional value of $\varphi_1 (P)$. 
Since  $\varphi_1 \!=\! \Psi$ in $\mathfrak{D}_3$, and 
$\varphi_1 \!=\! \Phi$ in $\mathfrak{D}_5$, 
we have for the exceptional values of $\varphi_1 (P)$
$$
\lambda^{'}_{2}, \,\lambda^{'}_{3}, \ldots, \,\lambda^{'}_{p}, \ldots,
$$
where the point set of $\varphi_1 (P) \!=\! \lambda^{'}_{p}$ is the same as  
the point set of $\varphi (P) \!=\! \lambda_p$.
Comparing  $\varphi_1 (P)$ with the original $\varphi (P)$, 
we easily see that 
 $\varphi_1 (P) \!=\! \varphi(P)$ in $\mathfrak{D}_{\mu_1}$, 
and $\varphi_1 (P) \!\geq\!  \varphi (P)$ in $\mathfrak{D}$.
Since $\varphi_1 \!\geq\!  \varphi$, 
$\varphi_1$ satisfies  Condition~$1^\circ$. 
This  $\varphi_1 (P)$ is a function satisfying almost the same property as
$\varphi (P)$. 
Although they differ only in the property of pseudoconvexity,
 the above operation
does not involve this property of $\varphi (P)$.
Therefore, in the same way as to produce $\varphi_1 (P)$
from $\varphi (P)$, 
we may construct $\varphi_2 (P)$ from $\varphi_1 (P)$.
We repeat this operation as far as the exceptional values remain,
and thus obtain
$$
\varphi_1 (P), \,\varphi_2 (P), \ldots, \,\varphi_p (P), \ldots .
$$
The part of properties of
$\varphi_p (P) \ (p \!>\! 1)$ which varies with $p$ is as follows: 
The exceptional values of $\varphi_p (P)$ are 
$$
\lambda^{(p)}_{p + 1}, \,\lambda^{(p)}_{p + 2}, \ldots, 
\,\lambda^{(p)}_{p + q}, \ldots ,
$$
where $\varphi_p (P) \!=\! \lambda^{(p)}_{p + q}$ and 
$\varphi (P) \!=\! \lambda_{p + q}$ are the same point set,
and in $\mathfrak{D}_{\mu_p}$,  
$\varphi_p (P) \!=\! \varphi_{p - 1} (P)$,
and in $\mathfrak{D}$, 
$\varphi_p (P) \!\geq\!  \varphi_{p - 1} (P)$ 
(note that in $\mathfrak{D}_5$, $\varphi_1 \!=\! \Phi$). 
We can thus choose such $\varphi_p (P)$.
Let $\varphi_0 (P)$ be the limit function of them, or the last function
in case the sequence is finite. 
Then $\varphi_0 (P)$ is clearly the required
function.  
\hfill C.Q.F.D.

\medskip
The function $\varphi_0 (P)$ thus obtained is in fact
a pseudoconvex function.\footnote{\,For this, 
the pseudoconvexity of $\mathfrak{D}_{\alpha_0} \!=\! \Delta$
suffices (Theorem 3 of the 9th Report).  Cf.\ \S9.} 

\bigskip
{\bf \S8.} At the beginning of the second
 Report\footnote{\,(Note by the translator.)
 This is the published second paper of the series
in J. Sci.\ Hiroshima Univ.\ Ser.\ A 7 (1937), 115--130.}
we explained the outer-convex ``H\"ulle'' with respect to
polynomials. 
We generalize it a bit more to supplement the fundamental lemma
of the previous section, but here we consider the (inner) convexity for
 convenience.

\medskip
{\bf Lemma 3.} {\em Let $\mathfrak{D}$ be a finitely sheeted
domain of holomorphy over $(x)$-space. 
Let  $E_0 \Subset \mathfrak{D}$ be an open subset.
Then, we have:\\
$1^\circ$  There exists a smallest open subset $H$
 among the open subsets of $\mathfrak{D}$, containing $E_0$,
which are convex with respect to the family of all holomorphic functions
in $\mathfrak{D}$,
and so $H \Subset \mathfrak{D}$. \\ 
$2^\circ$ There is no 
locally defined hypersurface $\sigma$ satisfying the following 
properties: 
$\sigma$ passes through a boundary point of 
$H$, but not through any  point of 
$H, \,E_0$ or the boundary of $E_0$, and the boundary points of $\sigma$ 
do not lie in $H$ nor on its boundary,
and $\sigma$ is defined in a form as follows: 
$$
\varphi (P) = 0, \quad P \in V, 
$$
where $V$ is a domain with $V \Subset \mathfrak{D}$, 
and $\varphi (P)$ is a holomorphic function in a neighborhood
of $V$ over 
$\mathfrak{D}$.
\footnote{\,Part $1^\circ$ above immediately follows from 
the existence theorem of  K--Konvexe H\"ulle due to
 H. Cartan--P.Thullen. 
Here, the original proof is based on a fundamental theorem of
simultaneous analytic continuation.
Cf.\ the paper of Cartan--Thullen mentioned above.
(See also the footnote of Theorem I.)}}

\medskip
{\it Proof.}~ $1^\circ$. 
We first show the existence of the H\"ulle $H$, for which we make 
some preparations.

Since $\mathfrak{D}$ is {\em finitely sheeted},
 a subset $\mathfrak{D}'$ of $\mathfrak{D}$ is bounded with respect to 
$\mathfrak{D}$ if and only if $\mathfrak{D}' \Subset \mathfrak{D}$.
Let $(\mathfrak{F})$ be the set of all holomorphic functions in
 $\mathfrak{D}$. 
Since $\mathfrak{D}$ is a domain of holomorphy, 
the First Theorem of {\em Cartan--Thullen}  
implies that $\mathfrak{D}$ is convex with respect to 
 $(\mathfrak{F})$. 
Therefore, regarding $\mathfrak{D} = \mathfrak{D}_0$ in
Lemma 1, we can construct an analytic polyhedron $\Delta$ of this lemma,
which is of the  form: 
$$
P \!\in\! R, \ | x_i| \!<\! r, \ | f_j (P)| \!<\! 1 \ \ (i \!=\! 1, 2, \ldots, 
n; j \!=\! 1, 2, \ldots, \nu). \leqno{\quad (\Delta)}
$$
Here, $f_j (P) \in (\mathfrak{F})$ 
and $R$ is an open subset of  
$\mathfrak{D}$ with $R \Supset \Delta$.
Further, note that for any given subset $E \Subset \mathfrak{D}$,
one may choose $\Delta \Supset E$. 

Let $\rho$ be an arbitrary positive number,
and let $d (P)$ denote the Euclidean boundary distance of
$\mathfrak{D}$. 
Let $\mathfrak{D}_\rho$ be the set of
all points $P \in \mathfrak{D}$ with $d (P) > \rho$.
(Here, $\rho$ is chosen so that $\mathfrak{D}_\rho$ is not empty.) 
If $\mathfrak{D}$ coincides with the finite $(x)$-space,
 then $\mathfrak{D}_\rho$ = $\mathfrak{D}$. 
By a parallel translation
$$
x'_{i} = x_i + a_i, \quad \sum | a_i|^2 \leq  \rho^2 \quad 
(i = 1, 2, \ldots, n),  \leqno{\quad (T)}
$$
we move a point $P$ of $\mathfrak{D}_\rho$ to $P'$ of $\mathfrak{D}$. 
If $P$ is given, 
$P'$ is uniquely determined.
For a function $f(P)$ of $(\mathfrak{F})$, we set
$$
F (P) = f (P').
$$
Then, $F (P)$ is a holomorphic function in $\mathfrak{D}_\rho$. 
Let $(T)$ be any of the parallel translation  within the restriction
mentioned above, and let $(\mathfrak{F}_\rho)$ be the set of all
functions $F(P)$ induced from  functions  $f (P)$ of 
$(\mathfrak{F})$. 

Let $A \Subset \mathfrak{D}$ be an open subset.  Assume that
$A$ is convex with respect to $(\mathfrak{F})$. 
\begin{itemize}
\item[]
{\em Let $A_0 \Subset A$ be an arbitrary open subset.
For a boundary point $M$ of $A$, there is a point $P_0$ arbitrarily
 close to $M$ such that there is at least one function
$f(P)$ of $(\mathfrak{F})$ with $| f (P_0)| >\max |f (A_0)|$.} 
\end{itemize}
We call this {\em Property $(\alpha)$} for a moment.
Conversely, we prove that if $A$ carries Property $(\alpha)$,
 {\em $A$ is convex with respect to $(\mathfrak{F})$}. 
Since  $A \Subset \mathfrak{D}$, 
an analytic polyhedron $\Delta$ above mentioned
is taken, so that $A \Subset \Delta$. 
Let $\rho$ be a sufficiently small positive number such that 
$\Delta \subset \mathfrak{D}_\rho$. 
Since $A$ satisfies Property $(\alpha)$,
 it is clear that $A$ is  convex with respect to $(\mathfrak{F}_\rho)$.
Now,  since every function of  $(\mathfrak{F}_\rho)$ 
is holomorphic in $\Delta$, 
it follows from Theorem 1 that it can be expanded to a series of
functions of $(\mathfrak{F})$,  
converging locally uniformly in $\Delta$.
Therefore, it is clear that $A$ is convex with respect to $(\mathfrak{F})$.

Now, let $A$ be an open subset of $\mathfrak{D}$, containing 
$E_0$ and  convex with respect to $(\mathfrak{F})$. 
Let $H$ be the subset of $\mathfrak{D}$ consisting of
all interior points of the intersection of all such $A$'s.

Since $E_0$ is open, $E_0 \subset H$. 
For $\Delta$ above, we may take 
$E = E_0$, and hence $H \Subset \mathfrak{D}$. 
It is clear that $H$ carries Property $(\alpha)$. 
Therefore, $H$ is convex with respect to $(\mathfrak{F})$. 
Thus, $H$ is the smallest open subset of $\mathfrak{D}$ which 
contains  $E_0$ and is convex with respect to $(\mathfrak{F})$, and
$H \Subset \mathfrak{D}$. 

\medskip
$2^\circ$.  
We assume the existence of a hypersurface  $\sigma$ with the properties
stated in the lemma.  
It suffices to deduce a contradiction.
Let $\varphi (P)$ be holomorphic in $V'$ such that 
$V \Subset V' \subset \mathfrak{D}$.
Let $d (P)$ denote the Euclidean boundary distance with respect to $V'$.
We choose a positive number $\rho$ such that
{$\min d (V) > \rho$
(the left-hand side of the inequality stands for the infimum of
$d(P)$ in $V$). }
Through the parallel translation
$$
x'_{i} = x_i + z_i, \quad \sum |z_i|^2 \leq  \rho^2 \qquad (i = 1, 2, \ldots, n),
$$
we move a point $P$ of $V$ to 
a point $P'$ of $V'$. Regarding $(z)$ as complex parameters, 
we set
$$
\psi (P, z) = \varphi (P'),
$$
and consider a family of hypersurface pieces,
$$
(\mathfrak{G}) ~:~\psi (P, z) = 0,\,P \in V.
$$ 
We take  $\rho$ small enough, so that 
the boundary of any hypersurface piece of $(\mathfrak{G})$ never intersects
$H$. 

Let $H_0$ be the set of all points $P$ of $H$ such that
$P$ does not belong to any hypersurface piece of $(\mathfrak{G})$.
Let $A_0$ be an open subset with $A_0 \Subset H_0$. 
As seen above, there is a minimal open subset $A$ of $\mathfrak{D}$,
containing $A_0$, which is convex with respect to $(\mathfrak{F})$.
Since $H$ is convex with respect to 
$(\mathfrak{F})$, similarly to the case of $H$ above,
 we have $A \Subset H$ by Lemma 1.
Now, we show that $A \subset H_0$. 

We describe a $2n$-dimensional ball
$S$ with radius $\rho$ and center at the origin in $(z)$-space. 
The open subset $(H, S) \ ((x) \in H, (z) \in S)$ in 
$(x, z)$-space is convex with respect to the set of all holomorphic functions
in the domain $(x) \in \mathfrak{D}$. 
Therefore by Theorem 2 there is a meromorphic function $G (P, Z)$
in $(H, S)$ such that it is congruent to 
$$
1 / \psi (P, z)
$$
in the intersection of $(H,S)$ and $(V,S)$, and it has no pole elsewhere.
(Theorem 2 is stated for finitely sheeted domains of holomorphy,
but in fact, it needs only the properties which are endowed with
$\mathfrak{D}_0$ in Lemma 1.) 

Suppose that $A$ is not contained in $H_0$. 
Then, $A$, which is an open set, contains a point outside $H_0$.
We may take a point $(z^0)$ in $S$ such that 
a point $P_0$ of $A$ lies on $\psi (P, z^0) = 0$.
With a complex variable $t$, we consider a function 

\medskip
\centerline{
$G (P, \,  t z^0)$.\footnote{\,
{(Note by the translator.)
In the manuscripts of Oka Library \cite{okap},
 References of Part I of the present article,
 this is misprinted as $G(P,\, t, z^0)$: It is confirmed
to see the 11th document, 1943, Catalogue of Dr. Kiyoshi Oka's
own handwriting manuscripts in the web-site of \cite{okap}.}}
}
\medskip

Then this is meromorphic when  $P$ is in $H$ and 
$t$ is in a neighborhood of the line segment $(0, 1)$, 
 has poles at $P = P_0, \,t = 1$, 
and $G (P, 0)$ has no pole in a neighborhood of $A$ (over $\mathfrak{D}$). 
As $t$ moves over the line segment $(0, 1)$ from $1$ to $0$, 
we denote by $t_0$  the last $t$ such that $G (P,\, t z^0)$ carries a pole in
$A$ or its boundary. 
Then, $G (P,\, t_0 z^0)$ has to carry a pole on the boundary 
of $A$ and to be holomorphic in $A$. 
et  $M$ be one of such poles. 
Let $P_1$ be a point of $A$, sufficiently close to $M$.
Since $A_0 \Subset H_0$ and $M$ is not a point of indeterminacy locus,
we have
$$
|G (P_1, \, t_0 z^0)| > \max |G (A_0, \, t_0 z^0)|.
$$

By Theorem 1, $G (P, \, t_0 z^0)$ is expanded to a series
of functions of $(\mathfrak{F})$, locally uniformly convergent in $A$:
This clearly contradicts the minimality of $A$. 
Thus, ``$A \subset H_0$'' holds. 

Since $A_0$ is an arbitrary open subset with 
$A_0 \Subset H_0$, the above consequence implies that
the open set $H_0$ satisfies  Property $(\alpha)$. 
Therefore, $H_0$ is convex with respect to $(\mathfrak{F})$; 
this conclusion holds no matter how $\rho$ is small.
Now, for sufficiently small $\rho$, $E_0 \subset H_0$: 
This again contradicts the minimality of $H$.
\hfill C.Q.F.D.

\bigskip
{\bf \S9.} The following two lemmata are easily deduced from Lemma 3.

\medskip
{\bf Lemma 4.} {\em Let $\Delta$ be a univalent domain of 
$(x)$-space which is  convex with respect to
polynomials, and let $\varphi (x)$
be a real-valued continuous function in a neighborhood of 
$\Delta$, satisfying Condition $2^\circ$ stated in Lemma II. 
If $\Delta_\alpha=\{x \in \Delta: \varphi (x) < \alpha\}$
for an arbitrarily given real number $\alpha$,
then $\Delta_\alpha$ is convex with respect to
polynomials, provided that it exists.}

\medskip
{\em Proof.}  It follows from Lemma 3 that there is a univalent minimal
open subset $H$ containing $\Delta_\alpha$,
which is convex with respect to polynomials.
Clearly,  $H \subset \Delta$. 
Therefore,  $\varphi (x)$ is defined in a neighborhood of $H$.
Let $\Bar{H}$ be the closure of $H$, and let $\beta$ be the
maximum value of $\varphi (x)$ on $\Bar{H}$. 
There are points of $\Bar{H}$ with $\varphi (x) = \beta$. 
Let $M$ be one of them. 
Since $\varphi (x)$ satisfies Condition~$2^\circ$, 
$M$ lies on the boundary of  $H$.
Furthermore, by the same property, there is a hypersurface
in a neighborhood of $M$, passing
through $M$ and no other point of $\Bar {H}$.
By Lemma 3, $M$ must be a boundary point of $\Delta_\alpha$. 
It follows that $\beta = \alpha$, and 
so $H = \Delta_\alpha$. 
Therefore, $\Delta_\alpha$  is convex with respect to
polynomials. \hfill{C.Q.F.D.}

\medskip
{\bf Lemma 5.} 
{\em Let $\varphi (P)$ be a real-valued continuous function
in a domain $\mathfrak{D}$ of $(x)$-space, satisfying Condition~$2^\circ$
in Lemma II. Let $\Delta$ be a domain of holomorphy such that 
$\Delta \Subset \mathfrak{D}$.
Put $\mathfrak{D}_\alpha=\{P \in \mathfrak{D}: \varphi (P) < \alpha\}$
for a real number $\alpha$.
If $\mathfrak{D}_\alpha \Subset \Delta$, 
then $\mathfrak{D}_\alpha$ is convex with respect to
all holomorphic functions in $\Delta$.
} 

\medskip
Since $\Delta \Subset \mathfrak{D}$, $\Delta$ is finitely sheeted.
Thus, $\Delta$ is a finitely sheeted domain of holomorphy,
and  $\mathfrak{D}_\alpha \Subset \Delta$.
Hence, Lemma 3 can be applied for $\Delta_\alpha$, and the rest
is exactly the same as above.

\medskip
We next state the theorems of {\em H. Cartan--P. Thullen} and 
{\em H. Behnke--K. Stein}: 

\medskip
{\bf The Second Theorem of H. Cartan--P. Thullen.}
{\em Let $\mathfrak{D}$ be a domain of $(x)$-space, and let 
$(\mathfrak{F})$ be the family of all holomorphic functions
in $\mathfrak{D}$. 
If the following two conditions are satisfied, then
 $\mathfrak{D}$ is a domain of holomorphy. 
\begin{itemize}
\item[$1^\circ$]
For an arbitrary set $\mathfrak{D}_0$ with
$\mathfrak{D}_0 \Subset \mathfrak{D}$,
there is an open set $\mathfrak{D}'$ with
 $\mathfrak{D}_0 \Subset \mathfrak{D}' \Subset \mathfrak{D}$
such that for every boundary point $M$ of $\mathfrak{D}'$ 
there is a function $f(P)$ of $(\mathfrak{F})$,
satisfying $|f (M)| > \max| f (\mathfrak{D}_0)|$.
\item[$2^\circ$]
For distinct two points $P_1, \,P_2$ of $\mathfrak{D}$, 
there is a function $f(P)$ of $(\mathfrak{F})$
with $f (P_1) \neq f (P_2)$.\footnote{\,The original authors stated
 this Second Theorem (also, the First Theorem)
in terms of $K$-convexity, but we stated it in the form above for
convenience: The proof is fully similar and direct.}
\end{itemize}
}

\medskip
{\bf Lemma of H. Behnke--K. Stein.} 
{\em Let $\mathfrak{D}$ be a domain of $(x)$-space, and let  
$$
\mathfrak{D}_1, \mathfrak{D}_2, \ldots, \mathfrak{D}_p, \ldots
$$
be a sequence of open subsets of  $\mathfrak{D}$ such that
 $\mathfrak{D}_p \Subset \mathfrak{D}_{p + 1}$  
and the limit is  $\mathfrak{D}$. 
We assume:
\begin{itemize}
\item[$1^\circ$]
Every $\mathfrak{D}_p$ is convex with respect to
the family $(\mathfrak{F} _{p + 1})$ of all holomorphic
functions in $\mathfrak{D}_{p + 1}$; 
\item[$2^\circ$]
For any two distinct points $P_1, P_2$ of $\mathfrak{D}_p$,
there is a function $f(P)$ in $(\mathfrak{F} _{p + 1})$ 
with $f(P_1) \neq f(P_2)$.
\end{itemize}
Then, $\mathfrak{D}_p$ has the same properties as
$1^\circ$ and $2^\circ$ above with respect to the family 
$(\mathfrak{F})$ of all holomorphic functions in $\mathfrak{D}$.%
\footnote{\,H. Behnke--K. Stein : 
Konvergente Folgen von Regularit\"atsbereichen 
und die Meromorphiekonvexit\"at, 1938 (Math.\ Annalen). \label{bk}}
}

\medskip
{\it Proof.} 
(Since $\mathfrak{D}_{p + 1}$ is a domain of holomorphy
by the Second Theorem of {\em Cartan--Thullen}), 
it follows from Theorem 1 
that every holomorphic function $\varphi (P)$ in $\mathfrak{D}_p$
is expanded to a series of functions of $(\mathfrak{F} _{p + 1})$,
locally uniformly convergent in  $\mathfrak{D}_p$.
This holds for  $p + 1, p + 2, \ldots$, as well,
and so $\varphi (P)$ may be similarly expanded to a series 
of functions of $(\mathfrak{F})$. 
Therefore, $\mathfrak{D}_p$ clearly has properties $1^\circ$ and $2^\circ$
with respect to $(\mathfrak{F})$.
\hfill C.Q.F.D.

\medskip
{\bf Theorem of H. Behnke--K. Stein.} 
{\em Let $\mathfrak{D}$ be a domain of $(x)$-space. 
Assume that for an arbitrary subset $\mathfrak{D}_0$ with
$\mathfrak{D}_0$ $\Subset \mathfrak{D}$, there is a domain
of holomorphy $\mathfrak{D}'$ with 
$\mathfrak{D}_0 \subset \mathfrak{D}' \Subset \mathfrak{D}$.
Then, $\mathfrak{D}$ is a domain of holomorphy.
\footnote{\,The same as \ref{bk}.} 
}

\medskip
{\it Proof.}
Since $\mathfrak{D}'$ is a domain of holomorphy,
it is pseudoconvex by {\em F.~Hartogs}.
Therefore it is inferred from Corollary 2 of Theorem 2
in the IX-th Report\footnote{\,(Note by the translator.)
This is the IX-th Report of the present series VII--XI, 1943.}
 that $\mathfrak{D}$ is pseudoconvex. 
Thus, there is a function $\varphi_0 (P)$ given in Lemma II for
$\mathfrak{D}$.  
By Lemma 5, $\mathfrak{D}_\alpha $
 ($\varphi_0 (P) < \alpha$, $P \in \mathfrak{D}$)
is convex with respect to all of holomorphic functions
in a domain of holomorphy $\mathfrak{D}'$ with
$\mathfrak{D}_\alpha \Subset \mathfrak{D}'$. 
Therefore, if $\alpha, \,\beta$ are arbitrary real numbers
with  $\alpha \!<\! \beta$, $\mathfrak{D}_\alpha$ satisfies the two conditions
stated in Lemma of {\em Behnke--Stein} 
with respect to all of holomorphic functions in
 $\mathfrak{D}_\beta$, and hence 
$\mathfrak{D}_\alpha$ satisfies the same with respect to
 all of holomorphic functions in $\mathfrak{D}$. 
Therefore by the Second Theorem of {\em Cartan--Thullen}, 
$\mathfrak{D}$ is  a domain of holomorphy. 
 \hfill C.Q.F.D.

We here generalize a bit more some parts of Lemmata 4 and 5.

\medskip 
{\bf Lemma 6.} {\em
Let $\mathfrak{D}$ be a finitely sheeted domain of holomorphy
over $(x)$-space, 
and let $\varphi (P)$ be a real-valued continuous function 
in $\mathfrak{D}$, satisfying Condition~$2^\circ$ in Lemma II.
If $\mathfrak{D}_\alpha =\{ P \in \mathfrak{D}: \varphi (P) < \alpha \}$
for  an arbitrarily given real number $\alpha$, 
then, every connected component of $\mathfrak{D}_\alpha$ is
a domain of holomorphy (provided that $\mathfrak{D}_\alpha$ is
not  empty).\footnote{\,In fact,
 $\mathfrak{D}_\alpha$ is convex for the family of all
holomorphic functions in  $\mathfrak{D}$. 
}
}

\medskip
{\it Proof.} Suppose that $\mathfrak{D}_\alpha$ exists.
Since $\mathfrak{D}$ is a domain of holomorphy, 
thanks to {\em F. Hartogs},  $\mathfrak{D}$ is pseudoconvex,
so that  there is a real-valued function $\psi (P)$ in
 $\mathfrak{D}$, stated in Lemma II.
Let $\beta$ be a real number with $\beta< \alpha$, 
and let  $\gamma$ be an arbitrary number. 
We consider an open set defined by
$$
P \in \mathfrak{D}, \quad \varphi (P) < \beta, \quad \psi (P) < \gamma.
 \leqno{\quad (\mathfrak{D}_{\beta \gamma})}
$$
Since $\mathfrak{D}$ is a finitely sheeted domain of holomorphy
and $\mathfrak{D}_{\beta \gamma} \Subset \mathfrak{D}$,
we can apply Lemma 3 with $E_0 = \mathfrak{D}_{\beta \gamma}$.
Hereafter, fully in the same way as the case of Lemma 4,
 we easily see that
$\mathfrak{D}_{\beta \gamma}$ is convex with respect to 
 all of holomorphic functions in $\mathfrak{D}$.
Therefore,  $\mathfrak{D}_{\beta \gamma} \subset \mathfrak{D}$,
so that by the Second Theorem of {\em Cartan--Thullen}, 
each connected component of $\mathfrak{D}_{\beta \gamma}$
is a domain of holomorphy. 
Note that $\mathfrak{D}_{\beta \gamma} \Subset \mathfrak{D}_\alpha$,
and $\mathfrak{D}_{\beta \gamma}$ can be chosen arbitrarily close
to $\mathfrak{D}_\alpha$. 
It follows from Theorem of {\it Behnke--Stein} that
each connected component of $\mathfrak{D}_\alpha$ 
is a domain  of holomorphy. 
\hbox{\quad} \hfill {C.Q.F.D.}\footnote{\,By F. Hartogs, 
domains of holomorphy are pseudoconvex, so that
we easily see the property of pseudoconvex domains by
Lemma 4 together with
 the theorems of the present section and those of the IX-th Report:
{\em Let $\varphi (x)$ be a pseudoconvex function in a
neighborhood of a $2 n$-dimensional ball $S$, and let 
$S_\alpha$ denote the sets of all points $x$ of $S$ with
$\varphi (x) < \alpha$  ($\alpha$ is an arbitrary real number).
Then, $S_\alpha$, if exists, is pseudoconvex.}
}

\bigskip
{\bf \S10.}   We return to our theme. In first, we claim that
a pseudoconvex domain is a domain of holomorphy.

We consider a {\em finitely sheeted} domain $\mathfrak{D}$
in $(x)$-space. We write
$$
x_1 = \xi + i \,\eta 
$$
with real and imaginary parts.
Let $a_1$ and $a_2$ be real numbers such that 
$$
a_2 < 0 < a_1,
$$ 
and denote by  $\mathfrak{D}_1$ the part of $\mathfrak{D}$
with $\xi \!<\! a_1$,
by  $\mathfrak{D}_2$ the part of $\mathfrak{D}$ with
$\xi \!>\! a_2$,
and set $\mathfrak{D}_3=\mathfrak{D}_1 \cap \mathfrak{D}_2$.
Assuming that the parts of $\mathfrak{D}$ with
 $\xi \!<\! a_2$ and $\xi \!>\! a_1$ are not empty,
we take points $Q_1, \,Q_2$ therein respectively.
Assume that {\em every connected component
of $\mathfrak{D}_1$ and $\mathfrak{D}_2$ is a domain of holomorphy}.
Then, necessarily so is $\mathfrak{D}_3$. 

Since a domain of holomorphy is pseudoconvex by {\em F. Hartogs},
$\mathfrak{D}$ is pseudoconvex. 
We may consider a real-valued function $\varphi_0 (P)$,
stated in Lemma~II for this $\mathfrak{D}$.
With a real number $\alpha$, 
we consider a subset $\mathfrak{D}_\alpha$ of $\mathfrak{D}$
such that $\varphi_0 (P) \!<\! \alpha$.
For a large $\alpha$, $\mathfrak{D}_\alpha$ contains the fixed points 
$Q_1$ and $Q_2$ in one connected component denoted by $A$.
It is noted that $A$ is bounded and finitely sheeted.
We denote respectively by $A_1, \,A_2, \,A_3$ 
the parts of $A$ with $\xi \!<\! a_1$, $\xi \!>\! a_2$,
and $a_2 \!<\! \xi \!<\! a_1$.
It follows from Lemma 6 that every connected component
of $A_1$, $A_2$ and $A_3$ is a domain of holomorphy.

We denote by $\Gamma$ the boundary of $A$ over $\xi \!=\! 0$.
Let $M$ be any point of $\Gamma$.
Then, there is a hypersurface piece $\sigma$ defined locally
in a neighborhood of $M$ and passing through $M$ such that
{$\sigma_0 \setminus\{M\}$ lies only in such
a part of a neighborhood of $\sigma$ in
$\mathfrak{D}$ that $\varphi_0(P) > \alpha$.}
Let $\psi (P) \!=\! 0$ ($\psi (P)$ is a holomorphic function)
be a defining equation of $\sigma$.
Choose $\beta$ with $\alpha \!<\! \beta$, sufficiently 
close to $\alpha$. 
Then $\sigma$ does not have the boundary point
in  $\mathfrak{D}_\beta \ (\varphi_0 (P) \!<\! \beta)$.
(Here, if necessary, we take out a neighborhood of the boundary
of $\sigma$.)
Let $B$ denote the part of $\mathfrak{D}_\beta$ with $a_2 \!<\! \xi \!<\! a_1$. 
Then, $B$ is finitely sheeted, and every connected component of $B$
is a domain of holomorphy.  Therefore, by Theorem 2
there is a function $G (P)$, meromorphic in $B$
with poles $1 / \psi (P)$ only on $\sigma$ and no other poles.
In $A_3$, $G(P)$ is holomorphic. 
For every point $M$ of $\Gamma$, 
there is  such a function $G (P)$ associated.  
Also, every connected component of $A_3$ is a domain of holomorphy
(cf.\ the method of the proof of Lemma 1).
Therefore, if positive $\delta_0$ and $\varepsilon_0$ are chosen
sufficiently small, by the standard arguments we easily deduce the
existence of holomorphic functions  
$f_j (P) \ (j \!=\! 1, 2, \ldots, \nu)$ in
 $A_3$ satisfying the following three conditions:
\begin{enumerate}
\item[{$1^\circ$}]  Let $A_0$ denote the set of all points of $A$
with $|\xi| \!<\! \delta_0, \, |f_j (P)| \!<\! 1$
$(j \!=\! 1, 2, \ldots, \nu)$. 
Then, $A_0 \Subset A$. 
\item[$2^\circ$]  Let $p$ be anyone of $1, 2, \ldots, \nu$. 
Then there is no point of $\mathfrak{D}_3$ with
$|f_p (P)| \!\geq\!  1 - \varepsilon_0$, 
lying over $|\xi - a_1| \!<\! \delta_0$, 
or over $|\xi - a_2| \!<\! \delta_0$. 
\item[$3^\circ$]  
The vector-valued function $[f_1 (P), f_2 (P), \ldots, f_\nu (P)]$
never takes the same vector-value for 
{mutually overlapped two points} of $A_0$.
\end{enumerate}

Further, letting $A_4$ be the set of points of $A_3$
satisfying $|f_j (P)| \!<\! 1 \ (j \!=\! 1, 2, \ldots, \nu)$,
we see that $A_4$ can be chosen arbitrarily close  to $A_3$. 
The union of $A_4$ and the part of $A$ satisfying $\xi \!\leq\!  a_2$ 
or $\xi \!\geq\!  a_1$ is an open set. 
Choose $f_j (P) \ (j \!=\! 1, 2, \ldots, \nu)$
so that $A_4$ is sufficiently close to $A_3$.
Then that open set contains the fixed points
 $Q_1$ and $Q_2$ in the same connected component,
which is denoted by $\Delta$. 
{\em The domain $\Delta$ satisfies the conditions given in \S4.}

If  $\alpha$ is chosen to be larger than a certain number $\alpha_0$, 
we may consider $A$ as a connected component of
$\mathfrak{D}_\alpha$, which contains $Q_1$ and $Q_2$. 
Choose $\alpha'$ with $\alpha_0 \!<\! \alpha' \!<\! \alpha$. 
In the same way as we associate $\alpha$ with $A$, 
we associate $\alpha'$ with $A'$. 
Needless to say,  $A' \Subset A$. 
Let  $A^{'}_{1}$ (resp.\ $A^{'}_{2}$) denote the part of $A'$ with 
$\xi \!<\! 0$ (resp.\ $\xi > 0$). 
Since $\Delta$ can be chosen arbitrarily close to $A$, 
we immediately obtain the following consequence from the result
 of the previous chapter:  
{\em For a given holomorphic function $\Phi (P)$ in the open
set, $P \!\in\! A$ with $|\xi| \!<\! \delta_0$ 
(here, $\delta_0$ can be arbitrarily small),
we can construct holomorphic function $\Phi_1 (P)$ (resp.\ $\Phi_2 (P)$)
in $A^{'}_{1}$ (resp.\  $A^{'}_{2}$), which is holomorphic in the part of $A'$
with $\xi \!=\! 0$,
such that $\Phi_1 (P) - \Phi_2 (P) = \Phi (P)$
holds there identically.} 

\medskip
Suppose that a pole  $(\wp)$ is given in $A$.
By Theorem 2 we may construct a meromorphic function $G_1 (P)$
in $A_1$ with pole  $(\wp)$. 
It is the same in $A_2$, and so  the meromorphic function is denoted by $G_2(P)$.
The difference $G_1 (P) - G_2 (P)$ is holomorphic 
in $A_3$. 
By the result above we see the following: {\em For a Cousin I Problem
given in $A$ we can solve it in $A'$.} 

\medskip
We come back to $A$: $A$ is a connected component of 
$\mathfrak{D}_\alpha \ (\alpha_0 \!<\! \alpha)$, 
containing $Q_1$ and $Q_2$. 
Let $M$ be any boundary point of $A$.
Let $(\gamma)$ be the polydisk described over 
$\mathfrak{D}$  with center $M$. 
For sufficiently small $(\gamma)$, 
there is a  hypersurface piece $\sigma$ defined
in $(\gamma)$, passing through $M$, which lies in
$\varphi_0 (P) \!>\! \alpha$ except for $M$.
Let $\sigma$ be defined by 
$$
\psi (P) = 0, \quad P \in (\gamma),
$$
where $\psi (P)$ is a holomorphic function in $(\gamma)$. 
If necessary, $(\gamma)$ is chosen a little smaller,
there is $\alpha''$ close to $\alpha$ with $\alpha \!<\! \alpha''$,
and the associated domain $A''$ contains no boundary point 
of $\sigma$. 
Therefore, by the arguments as above,
choosing  $\alpha''$ even closer to $\alpha$, 
we may obtain a meromorphic function 
 $G (P)$ in $A''$ such that it has poles 
$1 / \psi (P)$ over $\sigma$, 
and has no other pole.
Here $M$ is an arbitrary boundary point of $A$.

We examine the two conditions of the Second Theorem of {\em Cartan--Thullen}
for $A$. 
Let  $(\mathfrak{F})$ denote the set of all holomorphic functions
in $A$. 
Clearly by what we have seen above, 
$1^\circ$  $A$ is convex with respect to  $(\mathfrak{F})$.   

Let $P_1, P_2$ be an arbitrary pair of
{mutually overlapped points} of $A$
and denote the common base point by $\underline{P}$.
We describe a half-line $\underline{L}$ with one end at
 $\underline{P}$ in $(x)$-space. 
We describe a half-line on $A$ 
starting from $P_1$ over $\underline{L}$.
Since $A$ is bounded, this half-line necessarily intersects the boundary of  $A$.
Let $M_1$ be such a point, and let $L_1$
be the line segment $(P_1, M_1)$. 
Similarly, we describe a half-line $L_2$ starting from $P_2$.
Suppose that the length of $L_1$ does not exceed that
of $L_2$. 
(Clearly, this assumption does not lose generality.)
We denote by $G_0(P)$ the function $G(P)$ associated with
 $M = M_1$; 
$G_0 (P)$ is holomorphic in  $A$, holomorphic at every
boundary point of $A$ except for 
$M_1$, and has a pole at $M_1$. 
Therefore, $G_0 (P)$ has to have different function elements
at $P_1$ and $P_2$. 
Thus we have  2$^\circ$: For any distinct two points of $A$, 
there is necessarily a function of $(\mathfrak{F})$ having different values at
those points. 

Thus, Conditions 1$^\circ$ and $2^\circ$ are satisfied,
and so by the Second Theorem of {\em Cartan--Thullen}, 
$A$ is a domain of holomorphy. 
Since $\mathfrak{D}$ is a finitely sheeted domain, and
$A$ can be chosen arbitrarily close to it, 
{\em  Theorem of Behnke--Stein} implies that 
{\em $\mathfrak{D}$  is a domain of holomorphy.}  

\medskip
Now, we assume that $\mathfrak{D}$ 
is a {\em pseudoconvex domain} in $(x)$-space. 
For this $\mathfrak{D}$ we take a function
$\varphi_0(P)$ given in Lemma II, and consider
$\mathfrak{D}_\alpha \ (\varphi_0 (P) \!<\! \alpha)$
with an arbitrary real number $\alpha$.
(Here we take $\alpha$ enough large, so that
$\mathfrak{D}_\alpha$ really exists.) 
As in the proof of Theorem 2 (cf.\ \S3 and the last Report, \S3),
we divide  $\mathfrak{D}_\alpha$ into small 
$2 n$-dimensional cubes (open sets) $(A)$; here however, 
$(A)$ are not necessarily of complete form.
After sufficiently fine division, it follows from Lemma 4 that
every $(A)$ 
(not mentioning the case of complete form, but also in
another case) is a univalent open set, convex with respect to
polynomials.
Therefore, by the Second Main Theorem of {\em Cartan--Thullen}
every connected component of them is a domain of holomorphy.
After taking the division sufficiently fine,
it is the same for $(B)$ ($(B)_0$ is a $2n$-dimensional cube with center 
$(A)_0$, consisting of $9^n$ number of $(A)$ and 
some parts of their boundaries, which may be not of complete form).
Hence, from the result obtained above we easily infer in the same way
as in the case of {\em Cousin} I Problem  that every connected component
of $\mathfrak{D}_\alpha$ is a domain of holomorphy.
Therefore, Theorem of {\em Behnke--Stein} implies $\mathfrak{D}$
being a domain of holomorphy.

\medskip
{\bf Theorem I.}  {\em
A finite pseudoconvex domain with no interior ramification point is
 a domain of holomorphy.}\footnote{\,To detour
 around the use of the First Theorem of Cartan--Thullen,
 it suffices just
to replace ``{\em domain of holomorphy}'' by
 ``{\em domain $\mathfrak{D}$ satisfying the following two conditions}'':
Condition $1^\circ$, 
with $(\mathfrak{F})$ denoting the set of all holomorphic functions
in $\mathfrak{D}$, 
$\mathfrak{D}$  is convex with respect to $(\mathfrak{F})$.
; $2^\circ$, for every pair of distinct points of $\mathfrak{D}$ 
there is a function in $(\mathfrak{F})$ having distinct values
at the two different points. 
Consequently, Theorem I and  the  First Theorem of Cartan--Thullen
are obtained simultaneously.}

\medskip
By this theorem, the problem to show a domain being of holomorphy
is reduced to show the pseudoconvexity of
 the domain.\footnote{\,Cf.\ Report VI,
 Introduction.  As an example we frequently encounter,
we consider 
a ``\"Uberlargerungsbereich'' over a pseudoconvex domain, which is
pseudoconvex, too. 
Therefore, for example, in the Second Theorem of Cartan--Thullen,
the second condition is unnecessary.}

\bigskip
{\bf \S11.}  We extend the definition of convexity (the last Report,
\S1) a little, and redefine it as follows:

\medskip
{\bf Definition.} {\em Let $\mathfrak{D}$ be a finite domain over 
$(x)$-space with no interior ramification point,
 and let $(\mathfrak{F})$ be a family of holomorphic functions 
in $\mathfrak{D}$. 
The domain $\mathfrak{D}$ is said to be convex with respect to 
$(\mathfrak{F})$ if for 
every subset $\mathfrak{D}_0 \Subset \mathfrak{D}$,
there is an open set $\mathfrak{D}'$ with
$\mathfrak{D}_0 \subset \mathfrak{D}' \subset \mathfrak{D}$,
 bounded with respect to $\mathfrak{D}$, and
satisfying that for an arbitrary point
$P \in \mathfrak{D}\setminus \mathfrak{D}'$
there is at least one function $f(P)$ of $(\mathfrak{F})$
with $|f (P_0)| > \max |f (\mathfrak{D}_0)|$.
In the case where $\mathfrak{D}$ consists of finite or infinite
number of disjoint domains satisfying the property above,
we use the same terminologies as defined.} 

\medskip
The convexity in the sense of this definition clearly implies
that of the former definition.
It is convenient to consider the following convexity as well:

\medskip
{\bf Definition.} 
{\em In the above setting, $\mathfrak{D}$ is said to be strictly convex 
with respect to $(\mathfrak{F})$ if 
for every subset $\mathfrak{D}_0 \Subset \mathfrak{D}$,
there is an open set $\mathfrak{D}'$ with
$\mathfrak{D}_0 \subset \mathfrak{D}' \Subset \mathfrak{D}$,
satisfying the condition mentioned above.}   

\medskip
The strict convexity clearly implies the convexity.
If $\mathfrak{D}$ is finitely sheeted, 
these two new notions of convexity agree with the former one.
When  $\mathfrak{D}$ is convex (resp.\ strictly convex) 
with respect to the family of all
holomorphic functions in $\mathfrak{D}$,
$\mathfrak{D}$ is simply said to be 
{\em holomorphically convex} 
(resp.\ {\em strictly holomorphically convex}).\footnote{\,H. Behnke
 and people of his school use ``convexity''
in the sense of ``strict convexity''.
(Cf.\ Behnke--Thullen's monograph, 
the first two papers of H. Behnke--K. Stein referred at the beginning
of \S1, 
in particular the second one.) 
Here, as mentioned once before, the notion of global convexity
with respect to a family of holomorphic functions was introduced
by H. Cartan. (Cf.\ H. Cartan's paper referred in the footnote
at the end of Report IV.)}

It has been a question since the last Report
 if a domain of holomorphy is strictly
 holomorphically convex.\footnote{\,Cf.\ its \S1.
We did not leave from univalent domains until the first
research project (from Report I to Report VI) was finished:
The reason was at this point.}
We study it, here.

\medskip
{\bf Lemma 7.}  {\em 
In Lemma II {\rm(\S7)},  $\mathfrak{D}_\alpha$ is convex with respect to
the family of all holomorphic functions in $\mathfrak{D}$.
} 

\medskip
{\em Proof.} Note that $\mathfrak{D}_\alpha$
is pseudoconvex 
(due to Lemma 4, the Second Theorem of {\em Cartan--Thullen} 
and {\em Hartogs}' Theorem). 
Therefore, $\mathfrak{D}_\alpha$ is a domain of holomorphy by
Theorem I. 
Hence, with a real number $\beta$ such that $\alpha < \beta$, 
$\mathfrak{D}_\alpha$ is convex with respect to the family
of all holomorphic functions 
in $\mathfrak{D}_\beta$ by Lemma 5. 
Therefore, it follows from Lemma of {Behnke--Stein} that 
$\mathfrak{D}_\alpha$ is convex with respect to the family of all 
holomorphic functions in $\mathfrak{D}$.
\hfill  C.Q.F.D.

\medskip
{\bf Theorem II.}  {\em A finite domain of holomorphy is
strictly holomorphically convex.}

\medskip
{\em Proof.}~ Let $\mathfrak{D}$ be a (finite)
domain of holomorphy over $(x)$-space. 
Let  $E \Subset \mathfrak{D}$ be an arbitrary subset.
We take $\mathfrak{D}_\alpha$ in Lemma II so that 
 $E \Subset \mathfrak{D}_\alpha$. 
By Lemma 7 above, $\mathfrak{D}_\alpha$ is convex with respect to
the family of all holomorphic functions in  $\mathfrak{D}$,
and then by Lemma 1, with regarding
$\mathfrak{D}_0 = \mathfrak{D}_\alpha$,
we can choose an analytic polyhedron $\Delta$ of the form 
$$
P \!\in\! R, \ |x_i| \!<\! r, \ |f_j (P)| \!<\! 1 \ \ (i \!=\! 1, 2,\ldots, n; j \!=\! 1,
2,\ldots, \nu), \leqno{\quad (\Delta)}
$$
such that $E \Subset \Delta$.
Here,  $f_j (P)$ are functions of $(\mathfrak{F})$, 
and $R$ is a certain open set such that 
$\Delta \Subset R \subset \mathfrak{D}$. 

Let $P_0 \in \mathfrak{D}\setminus\Delta$ be any point.
It suffices to show that for this $P_0$ 
there is a function $f(P)$ of $(\mathfrak{F})$ 
with $|f (P_0)| > \max |f (E)|$. 
We take $\Delta'$ with the same property as $\Delta$ such that 
$\Delta \Subset \Delta'$ and $P_0 \in \Delta'$.
Let $\Delta'$ be of the form: 
$$
P \!\in\! R', \ |x_i| \!<\! r', \ |F_k (P)| \!<\! 1 \ \ (i \!=\! 1, 2,\ldots, n; k \!=\! 1,
2,\ldots, \mu). \leqno{\quad (\Delta')}
$$
Here, we choose $r'$ so that  $r \!\leq\!  r'$. 
From $\Delta$ and $\Delta'$ we form 
$$
P \!\in\! R', \quad |x_i| \!<\! r, \quad |f_j (P)| \!<\! 1, \quad |F_k (P)| \!<\! 1
\qquad \vspace{-2mm} 
\leqno{\quad (\Delta'')}
$$
$$
\qquad \qquad \qquad (i \!=\! 1, 2, \ldots, n;
 j \!=\! 1, 2, \ldots, \nu; k \!=\! 1, 2, \ldots, \mu)  .
$$
Clearly, $\Delta$ is one or a union of several connected components
of $\Delta''$. 
If $P_0$ does not belong to $\Delta''$, 
there exists necessarily a function with required property
among  $x_i, \,f_j (P)$. 
If $P_0$ belongs to  $\Delta''$, 
we consider a function in $\Delta''$ 
such that it is $0$ in $\Delta$, and $1$,
elsewhere. 
Then this function is holomorphic in $\Delta''$,
and so by Theorem 1 it is expanded to a series of functions 
of  $(\mathfrak{F})$, locally uniformly convergent in $\Delta''$.
Therefore, there is such a required function in this case, too.\\
 \hfill  C.Q.F.D.

\medskip
{\bf Corollary.}  {\em Let $\mathfrak{D}$ be a
 finite domain of holomorphy over $(x)$-space, and let 
$\mathfrak{D}_0$  be an open subset of $\mathfrak{D}$, convex 
with respect to the family of all holomorphic functions in
$\mathfrak{D}$.
Then, $\mathfrak{D}_0$ is strictly convex with respect
 to $(\mathfrak{F})$. 
}

\medskip
{\em Proof.}~ Since $\mathfrak{D}_0$ is convex with respect to
 $(\mathfrak{F})$, for any subset $E \Subset \mathfrak{D}_0$,
there is an open set $\mathfrak{D}'$ in  $\mathfrak{D}_0$ 
such that $E \subset \mathfrak{D}' \subset \mathfrak{D}_0$,
$\mathfrak{D}'$ is bounded with respect to 
$\mathfrak{D}_0$, 
and $\mathfrak{D}'$ satisfies the condition stated in the definition
of ``convexity''.
On the other hand, the above Theorem II implies the
existence of  an open set $\mathfrak{D}''$ in $\mathfrak{D}$
with $E \subset \mathfrak{D}'' \Subset \mathfrak{D}$,
which satisfies the same condition with respect to
$\mathfrak{D}$, and hence naturally with respect to
 $\mathfrak{D}_0$. 
We consider
 $\mathfrak{D}' \cap \mathfrak{D}'' = \mathfrak{D}_1$.
Then, $E \subset \mathfrak{D}_1 \subset \mathfrak{D}_0$
and satisfies this condition.
Now,  $\mathfrak{D}''$ is  finitely sheeted 
and $\mathfrak{D}'$ is bounded with respect to 
$\mathfrak{D}_0$, 
so that $\mathfrak{D}_1 \Subset \mathfrak{D}_0$.
Therefore, $\mathfrak{D}_0$ is strictly convex with respect to $(\mathfrak{F})$.
\hfill  C.Q.F.D.

\medskip
From Theorem 1 and this corollary,
we obtain the following consequence:

\medskip
{\bf Theorem III.} {\em
Let $\mathfrak{D}$ be a finite domain of holomorphy over  $(x)$-space, 
and let $\mathfrak{D}_0$ be an open subset of $\mathfrak{D}$,
which is convex
with respect to the family  $(\mathfrak{F})$ of all holomorphic
functions in $\mathfrak{D}$. 
Then, every holomorphic function in $\mathfrak{D}_0$
is expanded to a series of functions of $(\mathfrak{F})$,
convergent locally uniformly in  $\mathfrak{D}_0$.
} 

\medskip
The following result is deduced from
Theorem 2 and Theorems II and III:

\medskip
{\bf Theorem IV.} {\em 
In a finite domain of holomorphy, Cousin I Problem
is always solvable. 
} 

\begin{flushright}
(End, \  Report XI, ~3.12.12)  
\end{flushright}

\begin{flushright}
(Translated by Junjiro Noguchi (Tokyo))
\end{flushright}


\begin{thebibliography}{99}
\setlength{\itemsep}{-3pt}
\bibitem{an}
A. Andreotti and R. Narasimhan, Oka's Heftungslemma and the Levi problem,
Trans.\ Amer.\ Math.\ Soc.\ {\bf 111} (1964), 345--366.
\bibitem{bt}
H. Behnke and P. Thullen,
Theorie der Funktionen mehrerer komplexer Ver\"anderlichen,
Ergebnisse der Mathematik und ihrer Grenzgebiete Bd.\ 3,
Springer-Verlag, Heidelberg, 1934.
\bibitem{ct}
H. Cartan und P. Thullen,
Regularit\"ats- und Konvergenzbereiche, Math.\ Ann.\
{\bf 106} (1932), 617--647.
\bibitem{co}
P. Cousin, Sur les fonctions de $n$ variables complexes,
 Acta Math.\ {\bf 19} (1895), no. 1, 1--61.
\bibitem{br}
H.J. Bremermann,
 \"Uber die \"Aquivalenz der pseudokonvexen Gebiete und der
Holomorphiegebiete im Raum von $n$ komplexen Ver\"anderlichen,
Math.\ Ann.\ {\bf 128} (1954), 63--91.
\bibitem{forn}
J.E. Forn{\ae}ss,
A counterexample for the Levi problem for branched Riemann domains
over $\C^n$, Math.\ Ann.\ {\bf 234} (1978), 275--277.
\bibitem{fuj}
R. Fujita,
Domaines sans point critique int\'erieur sur l'espace projectif
complexe,
J. Math.\ Soc.\ Jpn.\ {\bf 15} (1963), 443--473.
\bibitem{grcas}
H. Grauert and R. Remmert,
Coherent Analytic Sheaves, Grundl.\ der Math.\ Wissen.\
vol.\ 265, Springer-Verlag, Berlin, 1984.
\bibitem{guro}
R.C. Gunning and H. Rossi,
Analytic Functions of Several Complex Variables, Prentice-Hall, 1965.
\bibitem{hi}
S. Hitotsumatsu, On Oka's Heftungs Theorem (Japanese), Sugaku
{\bf 1} (4) (1949), 304--307, Math.\ Soc.\ Jpn.
\bibitem{ho2}
L. H\"ormander,
Introduction to Complex Analysis in Several Variables,
First Edition 1966, Third Edition, North-Holland, 1989.
\bibitem{li}
I. Lieb, Das Levische Problem, Bonn.\ Math.\ Schr.\ No.\ 387 (2007), 1--34;
Translation into Frence, 
Le probl\`eme de Levi, Gaz.\ Math.\ Soc.\ Math.\ Fr.\
{\bf 115} (2008), 9--34.
\bibitem{nar}
R. Narasimhan, The Levi problem in the theory of several complex
	variables,
ICM Stockholm (1962).
\bibitem{nis96}
T. Nishino,
Function Theory in Several Complex Variables (in Japanese),
The University of Tokyo Press, Tokyo, 1996;
Translation into English by
N. Levenberg and H. Yamaguchi, Amer.\ Math.\ Soc.\,
Providence, R.I.,  2001.
\bibitem{nog16}
J. Noguchi,
Analytic Function Theory of Several Variables -- Elements of Oka's
 Coherence, Springer, Singapore, 2016;
translated from 
 Analytic Function Theory of Several Variables
(in Japanese), Asakura-Shoten, Tokyo, 2013.
\bibitem{nog19b}
J. Noguchi,
Analytic Function Theory of Several Variables (in Japanese,
 Tahensu Kaiseki Kansuron) -- Oka's Coherence for Undergraduates,
Second Edition,  Asakura-Shoten, Tokyo, 2019.
\bibitem{nog17}
J. Noguchi,
Inverse of Abelian  integrals and ramified Riemann domains,
 Math.\ Ann.\ {\bf 367} No.\ 1 (2017), 229--249.
\bibitem{nog18}
J. Noguchi, A weak coherence theorem and remarks to the Oka theory,
Kodai Math.\ J. {\bf 42} (2019), 566--586.
\bibitem{nog21}
J. Noguchi,  A New Introduction to Oka Theory---Basics of
Function Theory of Several Variables (in Japanese),
Shokabo, Tokyo (in press, 2021).
\bibitem{nog19}
J. Noguchi,
A brief chronicle of the Levi (Hartogs' Inverse) Problem, Coherence and
     an open problem, Manuscript, 2018,
Notices Intern.\ Cong.\ Chin.\ Math.\ {\bf 7} No.~2 (2019), 19--24.
\bibitem{nor}
F. Norguet,  Sur les domains d'holomorphie des fonctions uniformes
de plusieurs variables complexes (Passage du local au global),
Bull.\ Soc.\ Math.\ France {\bf 82} (1954), 137--159.
\bibitem{oka}
K. Oka:
\begin{itemize}
\item
Sur les fonctions analytiques de plusieurs variables -- I
 Domaines convexes par rapport aux fonctions rationnelles, 
  J.\ Sci.\ Hiroshima Univ.\ Ser.\ A {\bf 6} (1936), 245-255
 {\small[Rec.\ 1 mai 1936]}\footnote{\,The most commonly
cited reference for Oka's work should be \cite{oka2}, but
all the records of the received dates of the papers were
erased there by unknown reason: Here, they are listed
 for the sake of convenience.
}. 
\item
Sur les fonctions analytiques de plusieurs variables -- II Domaines d'holomorphie,
 J.\ Sci.\ Hiroshima Univ.\ Ser.\ A {\bf 7} (1937), 115-130
{\small[Rec.\ 10 d\'ec.\ 1936]}.
\item
Sur les fonctions analytiques de plusieurs variables -- III
 Deuxi\`eme probl\`eme de Cousin,
 J.\ Sci.\ Hiroshima Univ.\ {\bf 9} (1939), 7-19
{ [Rec.\ 20 jan.\ 1938]}.
\item
Sur les fonctions analytiques de plusieurs variables -- IV
 Domaines d'holomorphie et domaines rationnellement convexes, 
   Jpn.\ J.\ Math.\ {\bf 17} (1941), 517-521
{\small[Rec.\ 27 mar.\ 1940]}.
\item
Sur les fonctions analytiques de plusieurs variables -- V 
 L'int\'egrale de Cauchy, 
   Jpn.\ J.\ Math.\ {\bf 17} (1941), 523-531
{\small[Rec.\ 27 mar.\ 1940]}.\ 
\item
 Sur les domaines pseudoconvexes,
  Proc.\ of the Imperial Academy,
 Tokyo {\bf 17} (1941) 7-10 {\small[Communicated 13 jan.\ 1941]}.
\item
Sur les fonctions analytiques de plusieurs variables -- VI 
Domaines pseudoconvexes,
 T\^ohoku Math.\ J.\ {\bf 49} (1942), 15-52
{\small[Rec.\ 25 oct.\ 1941]}.
\item
Sur les fonctions analytiques de plusieurs variables -- VII 
 Sur quelques notions arithm\'etiques, 
 Bull.\ Soc.\ Math.\ France {\bf 78} (1950), 1-27
{\small[Rec.\ 15 oct.\ 1948]}.
\item
Note sur les fonctions analytiques de plusieurs variables,
  K{\= o}dai Math.\ Sem.\ Rep.\ (1949).\ no.\ {\bf 5-6}, 15--18
{\small[Rec.\ 19 d\'ec.\ 1949]}.
\item
Sur les fonctions analytiques de plusieurs variables -- VIII 
Lemme fondamental, 
 J.\ Math.\ Soc.\ Japan {\bf 3} (1951) No.\ 1, 204-214; 
 No.\ 2, 259-278 
{\small[Rec.\ 15 mar.\ 1951]}.\ 
\item
Sur les fonctions analytiques de plusieurs variables -- IX 
Domaines finis sans point critique int\'erieur, 
  Jpn.\ J.\ Math.\ {\bf 23} (1953), 97-155
{\small[Rec.\ 20 oct.\ 1953]}.
\item
Sur les fonctions analytiques de plusieurs variables -- X 
Une mode nouvelle engendrant les domaines pseudoconvexes, 
 Jpn.\ J.\ Math.\ {\bf 32} (1962), 1-12
{\small[Rec.\ 20 sep.\ 1962]}.
\end{itemize}
\bibitem{okapp}
K. Oka,  Appendice---Sur les formes objectives et les contenus subjectifs dans
 les sciences math\'ematiques; Propos post\'erieur---Pourquoi
 le pr\'esent m\'emoire est publi\'e de nouveau, 1953:
URL ``http://www.ms.u-tokyo.ac.jp/{\~{ }}noguchi/oka/''.
\bibitem{oka1}
K. Oka, Sur les fonctions analytiques de plusieurs variables,
Iwanami Shoten, Tokyo, 1961.
\bibitem{oka2}
K. Oka, Collected Works, Translated by
R. Narasimhan, Ed.\ R. Remmert, Springer-Verlag,
Berlin-Heidelberg-New York-Tokyo, 1984.
\bibitem{okap}
K. Oka,
Posthumous Papers of Kiyoshi Oka,
Eds.\ T. Nishino and A. Takeuchi, Kyoto, 1980--1983:
Oka Kiyoshi Collection, Library of Nara Women's University,
  URL ``http://www.lib.nara-wu.ac.jp/oka/index\_eng.html''.
\bibitem{tak}
A. Takeuchi,
Domaines pseudoconvexes infinis et la m\'etrique riemannienne
dans un espace projectif,
J. Math.\ Soc.\ Jpn.\ {\bf 16} (1964), 159--181.
\end{thebibliography}
\end{document}